\documentclass[12pt,reqno]{amsart}

\usepackage{latexsym}
\usepackage{graphicx}

\newcommand{\szego}{Szeg\"o }

\newcommand{\N}{{\mathbb N}}
\newcommand{\R}{{\mathbb R}}
\newcommand{\C}{{\mathbb C}}

\newcommand{\vol}{{\operatorname{Vol}}}

\newcommand{\supp}{{\operatorname{Supp\,}}}
\renewcommand{\phi}{\varphi}

\newcommand{\ispa}[1]{\langle \,#1 \,\rangle } 
\newcommand{\mf}{\mathfrak} 
\newcommand{\mb}{\mathbb} 
\newcommand{\ol}{\overline} 
\newcommand{\lspan}{{\rm span}} 
\newcommand{\sgn}{{\rm sgn}} 

\newcommand{\pcal}{\mathcal{P}}

\newtheorem{theo}{{\sc Theorem}}[section]
\newtheorem{maintheo}{{\sc Theorem}}
\newtheorem{cor}[maintheo]{{\sc Corollary}}

\newtheorem{lem}[theo]{{\sc Lemma}}
\newtheorem{prop}[theo]{{\sc Proposition}}
\newenvironment{example}{\medskip\noindent{\it Example:\/} }{\medskip}
\newenvironment{rem}{\medskip\noindent{\it Remarks:\/} }{\medskip}

\title[Lattice path combinatorics]
{Lattice path combinatorics and asymptotics of  multiplicities  of
weights in tensor powers}

\author{Tatsuya Tate}
\address{Department of Mathematics, Keio University, Keio University
3-14-1 Hiyoshi Kohoku-ku, Yokohama, 223--8522 Japan}
\email{tate@math.keio.ac.jp}

\author{Steve Zelditch }
\address{Department of Mathematics, Johns Hopkins University, Baltimore, MD
21218, USA}
\email{zelditch@math.jhu.edu}

\thanks{Research partially supported by JSPS (first author).}
\thanks{Research partially supported by NSF grants DMS-0071358 (second author).}

\date{\today}

\begin{document}

\begin{abstract} We give  asymptotic formulas for the multiplicities of weights and
irreducible summands in high-tensor powers $V_{\lambda}^{\otimes
N}$ of an irreducible representation $V_{\lambda}$ of a compact
connected Lie group $G$. The weights are allowed to depend on $N$,
and we obtain several regimes of pointwise  asymptotics, ranging
from a
 central limit region to  a large deviations region. We use a
 complex steepest descent method that applies to general
 asymptotic counting problems for lattice paths with steps in a convex
 polytope.
\end{abstract}

\maketitle

\tableofcontents

\section*{Introduction}

This article is concerned with the interplay between combinatorics
of lattice paths with steps in a convex polytope and asymptotics
of weight multiplicities (and multiplicities of irreducible
representations)  in high tensor powers $V_{\lambda}^{\otimes N}$
of irreducible representations $V_{\lambda}$ of a compact
connected Lie group $G$. Our main results give asymptotic formulae
for
\begin{itemize}
\item multiplicities $m_{N}(\lambda;\nu)$ of  weights $\nu$
in  $V_{\lambda}^{\otimes N}$;

\item  multiplicities  $a_{N}(\lambda;\nu)$ of irreducible representations $V_{\nu}$ with
highest weight $\nu$ in $V_{\lambda}^{\otimes N}$.

\item multiplicities of lattice paths with steps in a convex
lattice polytope $P$ from $0$ to an $N$-dependent lattice point
$\alpha \in NP.$

\end{itemize}
Asymptotic analysis of multiplicities in high tensor powers  are
of interest because  the known formulae for  multiplicities of
weights and irreducibles in tensor products (Steinberg's formula,
Racah formula, Littlewood-Richardson rule and others \cite{FH,
BD}) rapidly become complicated as the number of factors
increases.

Our analysis of multiplicities is based on the simple and
well-know fact \cite{S} that the multiplicities of lattice paths
can be obtained as Fourier coefficients of powers $k(w)^N$ of a
complex exponential sum of the form
\begin{equation}
k(w) = \sum_{\beta \in P}c(\beta)e^{\ispa{\beta, w}}, \quad w \in
\C^n
\end{equation}
with positive coefficients $c(\beta)$, where $P$ is a convex
lattice polytope. One can obtain the precise asymptotics of the
Fourier coefficients of $k(w)^N$ by  a complex stationary phase
(or steepest descent) argument. It is necessary to deform the
contour of the Fourier integral to pick up the relevant complex
critical points and to study the geometry of the complexified
phase, which is closely related to the moment map for a toric
variety. In fact, it was the analysis of this latter problem in
\cite{TSZ1, SZ2} which led to the present article.

When $P$ is the convex hull of a Weyl orbit of the weight
$\lambda$, the Fourier coefficients are weights of
$V_{\lambda}^{\otimes N}$. When $P=p\Sigma$ with the simplex
$\Sigma$ and a positive integer $p$, and $c(\beta)={p \choose
\beta}$ ($|\beta| \leq p$), then the Fourier coefficients are of
course multinomial coefficients of the form ${Np \choose \gamma}$
with $|\gamma| \leq Np$. Thus, lattice path multiplicities in
general behave much like multinomial coefficients, whose
asymptotics (obtained form Stirling's formula) have been studied
since Boltzmann in probability theory and statistical mechanics
(cf. \cite{E, F}). In view of the rather basic nature of the
lattice path counting problem and its applications, it might seem
surprising that a pointwise asymptotic analysis has not  been
carried out before (at least, to our knowledge). The closest prior
result appears to be Biane's central limit asymptotics for
multiplicities of irreducibles in tensor products \cite{B}, which
does not make use of the connection to lattice path counting.

To state our results, we  need some notation. We  fix a maximal
torus $T \subset G$ and  denote by $\mf{g}$ and $\mf{t}$ the
corresponding Lie algebras. Their duals are denoted by
$\mf{g}^{*}$ and $\mf{t}^{*}$. We fix an open Weyl chamber $C$ in
$\mf{t}^{*}$, and denote the set of dominant weights by $I^* \cap
\overline{C}$ where $I^*$ is the lattice of integral forms in
$\mf{t}^{*}$. For $\lambda \in I^* \cap \overline{C}$, we denote
by  $V_{\lambda}$ the irreducible representation of $G$ with the
highest weight $\lambda$, and denote its character by
$\chi_{V_{\lambda}}$ or more simply by $\chi_{\lambda}$.   We
further denote by $Q(\lambda) \subset \mf{t}^{*}$ the convex hull
of the orbit of $\lambda$ under the action of the Weyl group $W$.
The multiplicity of a weight $\mu$ in $V_{\lambda}$ is denoted by
$m_{1}(\lambda;\mu)$. We set $ M_{\lambda}=\{\mu
\,;\,m_{1}(\lambda;\mu) \neq 0\} \subset Q(\lambda).$

It is well known that the weights (and highest weights of
irreducibles) occurring in $V_{\lambda}^{\otimes N}$ all lie
within $Q(N \lambda).$  Our aim is to obtain pointwise asymptotic
formulae for the multiplicities for all possible weights. As will
be seen, the asymptotics fall into several regimes. We begin with
some simple  results on the bulk properties of  weight asymptotics
and progress to our main results giving individual asymptotic
formulae.

The simplest problem is to  determine the asymptotic distribution
of multiplicities of weights in $V_{\lambda}^{\otimes N}$. Let us
define a probability measure on $Q( \lambda)$ as follows:
\begin{equation}
d m_{\lambda, N} := \frac{1}{\dim V_{\lambda}^{\otimes N}}
\sum_{\nu \in Q(N \lambda)} m_N(\lambda, \nu) \delta_{N^{-1} \nu}.
\end{equation}
This measure charges each possible weight $\nu$  of
$V_{\lambda}^{\otimes N}$ with its relative multiplicity
$\frac{m_N(\lambda, \nu)}{\dim V_{\lambda}^{\otimes N}}$ and then
dilates the weight back to $Q(\lambda)$. As $N \to \infty$, the
dilated weights become denser in $Q(\lambda)$ and we may ask how
they become distributed. In particular,  which are the most
probable weights?

\begin{maintheo}  Assume that $\lambda$ is a dominant weight in the open Weyl chamber.
Then, we have
\[
m_{\lambda,N} \to \delta_{Q^{*}(\lambda)}
\]
weakly as $N \to \infty$, where $\delta_{Q^{*}(\lambda)}$ is the
Dirac measure at the (Euclidean)  center of mass $Q^{*}(\lambda)$
of the polytope $Q(\lambda)$ given by
\begin{equation}
Q^{*}(\lambda) =\frac{1}{\dim V_{\lambda}}
\sum_{\nu \in M_{\lambda}}m_{1}(\lambda;\nu)\nu.
\label{COMrep}
\end{equation}
\label{SIMPLE}
\end{maintheo}

This is an elementary result because
\begin{equation} \chi_{V_{\lambda}^{\otimes N}} =
\chi_{V_{\lambda}}^N \implies  dm_{\lambda, N} = D_{\frac{1}{N}}
dm_{\lambda}
* \cdots * dm_{\lambda},
\end{equation}
where $dm_{\lambda} = dm_{\lambda, 1}$ and where $D_{\frac{1}{N}}$
is the dilation operator by $\frac{1}{N}$ on the
dual Cartan subalgebra $\mf{t}^{*}$.
Hence, the  sequence of measures $\{dm_{\lambda, N} \}$ satisfies
the central limit theorem and  the (Laplace) large deviations
principle. Let us recall the definitions: Let $m_{N}$
($N=1,2,\ldots$) be a sequence of probability measures on a closed
set $E \subset \mb{R}^{n}$. Let $I:E \to [0,\infty]$ be a lower
semicontinuous function. Then, the sequence $m_{N}$ is said to
satisfy the {\it large deviation principle with the rate function}
$I$ (and with the speed $N$) if the following conditions are
satisfied:
\begin{itemize}
\item[(1)] The level set $I^{-1}[0,c]$ is compact for every $c \in \mb{R}$.
\item[(2)] For each closed set $F$ in $E$,
\[
\limsup_{N\to \infty}\frac{1}{N}\log m_{N}(F) \leq -\inf_{x \in
F}I(x).
\]
\item[(3)] For each open set $U$ in $E$,
\[
\liminf_{N \to \infty}\frac{1}{N}\log m_{N}(U) \geq -\inf_{x \in
U}I(x).
\]
\end{itemize}

The following is a consequence of Cram\'{e}r's theorem (\cite{DZ},
Theorems 2.2.3, 2.2.30):

\begin{maintheo} \label{CRAMER}
Assume that $G$ is semisimple.
Then, the sequence $\{dm_{\lambda, N}\}$ of measures on
$Q(\lambda)$ satisfies a large deviations principle with speed $N$
and rate function:
\begin{equation}
I_{\lambda}(x) = \sup_{\tau \in \mf{t}}
\left\{ \ispa{\tau, x} -
\log \left(
\frac{\chi_{\lambda}(\tau/(2\pi i))}{\dim V_{\lambda}}
\right)
\right \}
,\quad x \in \mf{t}^{*},
\label{RATEG}
\end{equation}
where $\chi_{\lambda}(\tau/(2\pi i))=\sum_{\nu \in M_{\lambda}}m_{1}(\lambda;\nu)e^{\ispa{\nu,\tau}}$
denotes the character of $V_{\lambda}$ extended on $\mf{t} \otimes \mb{C}$.
\end{maintheo}

The assumption that $G$ is semisimple is not necessary. However, in general case,
the definition of the rate function is slightly modified. See Section \ref{GREP} for details.

Before stating our more refined results, we consider the analogous
statements for multiplicities of irreducibles.  We  define the
analogous measures weighting $\mu \in Q(N \lambda)$ by the
multiplicity of the irreducible representation $V_{\mu}$ in
$V_{\lambda}^{\otimes N}$, defined by
\begin{equation}
dM_{\lambda, N} := \frac{1}{B_{N}(\lambda)}
\sum_{\nu \in Q(N \lambda)} a_N(\lambda, \nu) \delta_{N^{-1} \nu},
\quad B_{N}(\lambda)=\sum_{\mu}a_{N}(\lambda;\nu).
\end{equation}
The measures $dM_{\lambda,N}$ are measures on the closed positive
Weyl chamber $\overline{C}$. They also satisfies the Laplace large
deviations principle, but the proof is not quite as simple as for
$dm_{\lambda, N}$.  The measures $dM_{\lambda, N}$ and
$dm_{\lambda, N}$ are related by an alternating  sum over the Weyl
group (see Proposition \ref{Mofirred} and Lemma \ref{ALTM}).

\begin{equation}
dM_{\lambda, N}(\mu) = \frac{(\dim
V_{\lambda})^{N}}{B_{N}(\lambda)} \sum_{w \in W}\sgn (w)
dm_{\lambda, N}(\mu +\rho -w\rho).
\end{equation}

We can thus deduce the upper-bound half (see (2)) in the large
deviation principle for the measure $dM_{\lambda,N}$ from that for
$dm_{\lambda,N}$. It follows from Theorem \ref{CRAMER} that:

\vspace{5pt}

\begin{cor} \label{CRAMERM}
Assume that $G$ is semisimple.
The sequence $\{dM_{\lambda, N}\}$ of measures on
$Q(\lambda)$ satisfies the upper-bound in a large deviations principle with speed $N$
and rate function $I_{\lambda}(x)$ given by \eqref{RATEG}.
\end{cor}

The lower bound will follow from our pointwise asymptotics. We
should note the large deviations principle with the rate function
\eqref{RATEG} has already been  proved by Duffield \cite{D} for
$dM_{\lambda, N}$ by a different method.

These results give  the bulk properties of the measures
$dm_{\lambda, N}, dM_{\lambda, N}$ in that they give the exponents
of the measures of $N$-independent closed/open sets. Our main
results give apparently optimal refinements, in which we  give
pointwise asymptotics for multiplicities of ($N$-dependent)
weights.  As mentioned above, they are based on the combinatorics
of lattice paths rather than on large deviations theory, which
does not seem capable of seeing the finer details of the
asymptotics.

 To introduce our
results, we recall one of the first and most basic results of a
similar kind, namely Boltzmann's analysis of the asymptotics of
multinomial coefficients (see \cite{E} for historical background
and the relation to the present problem):
$$\left\{ \begin{array}{l} m_N: \{(k_1, \dots, k_m) \in \N^m: k_1 + \cdots + k_m \leq N \} \to \R^+, \\ \\
\;\;\;m_N(k_1, \dots, k_n)
 = {N \choose k_1, \dots, k_m} = \frac{N!}{k_1! \cdots k_m!}. \end{array} \right.$$

 Let us consider the case $m = 1$ of binomial
coefficients. It is easy to see that the binomial coefficient
$b_N(k) = {N \choose k}$ peaks at the center $k = \frac{N}{2}$ and by Stirling's
 formula $r! \sim \sqrt{2 \pi}  r^{r + \frac{1}{2}} e^{-r}$, $b_N(\frac{N}{2}) \sim N^{-1/2} 2^N$.
We measure distance from the center by   $d_{N}(k) = k -
\frac{N}{2}.$ We then have (see \cite{F}, Chapter 7 for the first
two lines ):
\[
b_N(k) \sim \left\{
\begin{array}{ll}
(CL)  \;\;\;  C \; N^{-1/2} 2^N \; e^{-\frac{2d_{N}(k)^2}{N}},
& \; \mbox{if}\; d_{N}(k) = o(N^{\frac{2}{3}})\\
& \\
(MD) \;\;\;  C\; N^{-1/2} 2^N \; e^{-\frac{2d_{N}(k)^2}{N} -
Nf(\frac{2d_{N}(k)}{\sqrt{N}})}\;\;, & \; \mbox{if}\; d_{N}(k) = o( N), \\ & \\
\phantom{(MD) \;\;\;  }\mbox{with}\;\;
f(x) = \sum_{n = 2}^{\infty}
\frac{x^{2n}}{(2n)(2n-1)}  &  \\ & \\
(SD) \;\;\; \frac{1}{\sqrt{2\pi Na(1 - a)}} \; a^{- a N} \; (1 - a)^{-
(1 - a) N}, & \;\; k \sim a N, \;\; a < 1; \\  & \\
(RE)\;\;\; C_0\; N^{k_0} , & k = k_0, \;  N - k_0.
\end{array} \right.
\]
We refer to the first region as the central limit region (CL),
where the asymptotics are normal (i.e. have the form $ N^{-1/2}
2^N  \phi (\frac{d_{N}(k)}{\sqrt{N}}),$ where $\phi$ is the
Gaussian). The exponential growth  is fixed at $\log 2$  as long
as $d_{N}(k) = O( \sqrt{N})$.   In the next region (MD) of
moderate deviations, the exponent is decreased  by the function
$f$. In the next regime (SD) of strong deviations, the growth
exponent is $a \log \frac{1}{a} + (1 - a) \log
\frac{1}{1 - a} < \log 2$. In the final boundary (RE) region
 of rare events, the exponent vanishes and  the growth rate is
algebraic.

In a somewhat similar way, multiplicities peak at weights near the
center of gravity $Q^*(\lambda)$  of $Q(N \lambda)$, have a common
exponential rate for weights  in a ball of radius $O(\sqrt{N})$
around the center of mass, and then the exponential rate declines
as the weight moves from a moderate to a strong deviations region
towards the boundary of $Q(N \lambda). $ At the boundary point  $N
\lambda$ of $Q(N \lambda)$, the multiplicity equals one.

\subsection{Statements of results on weight multiplicities}

To state our results precisely, we will need further notation.
Let $X^{*} \subset \mf{t}^{*}$ denote the subspace spanned by the
simple roots, and let $X=(X^{*})^{*}$ be its dual space. Using an
inner product which is invariant under the action of the Weyl
group, the space $X$ is identified with the subspace of $\mf{t}$
spanned by the inverse roots. As is shown in Section \ref{GREP},
the polytope $Q(\lambda)-\lambda$ is contained in $X^{*}$. In the
following, the interior of $Q(\lambda)$ means the interior of
$Q(\lambda)$ in the affine subspace $X^{*} +\lambda$. Let $\rho$
denote half the sum of the positive roots. Let $L^{*}$ be the
lattice of weights in $X^{*}$. Since all the roots is in $L^{*}$,
the lattice $L^{*}$ is of maximal rank in $X^{*}$.
Let $\Lambda^{*}$ be the root lattice in $X^{*}$, i.e., $\Lambda^{*}$
is the linear span of all the roots over $\mb{Z}$, which satisfies
$\Lambda^{*} \subset L^{*}$.
The both lattices $\Lambda^{*}$ and $L^{*}$ are of maximal rank.
Their duals are denoted by $\Lambda$ and $L$ respectively.
Then we have $L \subset \Lambda$, and hence the quotient
$\Pi(G) := \Lambda/L$ is a finite abelian group.

\subsubsection{Central limit region}

Our first result concerns the `central limit region' of weights
which are within a ball of radius $O(\sqrt{N})$ around the center
of mass in Theorem \ref{SIMPLE}. For the sake of simplicity we
will assume that $G$ is semisimple. In this case, we have $X^{*}=\mf{t}^{*}$,
and we can use the (negative) Killing form for the inner product
invariant under the action of the Weyl group.

\begin{maintheo} \label{CLR}
Assume that $G$ is semisimple. Fix a dominant weight $\lambda$ in the open Weyl chamber $C$.
Let $\nu$ be a weight such that $|\nu|=O(N^{1/2})$.
Assume that $m_{N}(\lambda;\nu) \neq 0$ for every
sufficiently large $N$. Then, we have
\begin{equation}
m_{N}(\lambda;\nu)= (2\pi N)^{-m/2} |\Pi(G)|(\dim
V_{\lambda})^{N} \left(
\frac{e^{-\ispa{A_{\lambda}^{-1}\nu,\nu}/(2N)}}{\sqrt{\det
A_{\lambda}}} +O(N^{-1/2}) \right),
\label{CLTmult}
\end{equation}
where $|\Pi(G)|$ is the order of the
finite group $\Pi(G)=\Lambda/L$, $m=\dim \mf{t}$ is the rank of $G$
and the positive definite
linear transform $A_{\lambda}:\mf{t} \to \mf{t}^{*}$ is given by
\begin{equation}
A_{\lambda}=\frac{1}{\dim V_{\lambda}}\sum_{\mu \in
M_{\lambda}}m_{1}(\lambda;\mu)\mu \otimes \mu. \label{repA}
\end{equation}
\label{repCLT}
\end{maintheo}

We note that in this regime, the exponent of growth of
multiplicities is the constant $\log \dim V_{\lambda}$. The
assumption that $m_{N}(\lambda;\nu) \neq 0$ for every sufficiently
large $N$ can be replaced by that $m_{N_{0}}(\lambda;\nu) \neq 0$
for some $N_{0}$ if $0$ is a weight in $V_{\lambda}$. In \S 2, we
prove a stronger result, Theorem \ref{MDgen}, which extends the
central limit regime to   weights  $\nu \in NQ(\lambda)$ of the
form
\begin{equation} \label{COG}
\nu =NQ^{*}(\lambda) +d_{N}(\nu),\quad |d_{N}(\nu)|=o(N^{s}).
\end{equation} with $0 \leq s \leq 2/3$. Here, as in the case of binomial
coefficients, $d_{N}(\nu)$ represents the distance to
 the center of gravity of
$Q(\lambda)$.

\subsubsection{Large deviations region}

We now consider the moderate and strong deviations regimes. As
suggested by the behavior of multinomial coefficients, the
exponent must decrease as we move away from the center of gravity
of $Q(N \lambda)$. A key role in the exponent correction will be
played by the map:
\begin{equation} \mu_{\lambda}: X \to Q(\lambda),\;\;
\mu_{\lambda}(x): =  \frac{1}{\sum_{\mu \in
M_{\lambda}}m_{1}(\lambda;\mu)e^{\ispa{\mu,x}}} \sum_{\mu \in
M_{\lambda}}m_{1}(\lambda;\mu)e^{\ispa{\mu, x}}\mu.
\label{Gmomexp}
\end{equation}
This map is a homeomorphism from $X$ to the interior of
$Q(\lambda)$ (see e.g. \cite{Fu}), and resembles the moment map of
a toric variety, restricted to the real torus in $(\C^*)^m$.
We define a function $\delta_{\lambda}$ on the interior $Q(\lambda)^{o}$
of the polytope $Q(\lambda)$ by
\begin{equation}
\delta_{\lambda}(x)=\log \left( \sum_{\mu \in
M_{\lambda}}m_{1}(\lambda;\mu)e^{\ispa{\mu - x
,\tau_{\lambda}(x)}} \right),
\label{Gconst0}
\end{equation}
where $\tau_{\lambda}=\mu_{\lambda}^{-1}:Q(\lambda)^{o} \to X$.
It is clear that $\delta_{\lambda}(\nu) >0$ for $\nu \in Q(\lambda)^{o} \cap M_{\lambda}$.
When $G$ is semisimple, the function $\delta_{\lambda}$ is related to
the rate function $I_{\lambda}$ given by \eqref{RATEG} by the formula
\begin{equation}
\label{RvsD}
\delta_{\lambda}(x)=\log (\dim V_{\lambda})-I_{\lambda}(x),
\quad x \in Q(\lambda)^{o}.
\end{equation}

For $\nu \in Q(\lambda)^{o}$, we further define the linear map
$A_{\lambda}^{0}(\nu):\mf{t} \to \mf{t}^{*}$ by
\begin{equation}
A_{\lambda}^{0}(\nu)=\sum_{\mu \in M_{\lambda}}
\frac{m_{1}(\lambda;\mu)e^{\ispa{\mu,\tau_{\lambda}(\nu)}}}
{\sum_{\mu' \in M_{\lambda}}m_{1}(\lambda;\mu')e^{\ispa{\mu',\tau_{\lambda}(\nu)}}}
\mu \otimes \mu -\nu \otimes \nu.
\label{Gmatrix00}
\end{equation}
In general, the linear transform $A_{\lambda}^{0}(\nu)$ defined above has a zero eigenvalue.
However, its restriction to the subspace $X$, which is denoted by
\begin{equation}
A_{\lambda}(\nu):=A_{\lambda}^{0}(\nu)|_{X},
\label{Gmatrix0}
\end{equation}
is shown to be positive definite as a linear map from $X \to X^{*}$.

First, we consider the `strong deviations' regime where  the weight
in question has the form $\nu = N \nu_0 + f$.

\begin{maintheo}
Let $\lambda \in C \cap I^{*}$ be a dominant weight,  and let
$\nu_0 \in M_{\lambda}$ be a weight of $V_{\lambda}$ which lies in
the interior $Q(\lambda)^{o}$ of the polytope $Q(\lambda)$. We fix
a weight $f$ in the root lattice $\Lambda^{*}$. Then, we have the
following asymptotic formula:
\[
m_{N}(\lambda;N\nu_0 +f)= (2\pi N)^{-m/2}
\frac{|\Pi(G)|e^{N\delta_{\lambda}(\nu_0)
-\ispa{f,\tau_{\lambda}(\nu_0)}}} {\sqrt{\det A_{\lambda}(\nu_0)}}
(1+O(N^{-1})),
\]
where $m$ is the number of the simple roots, $|\Pi(G)|$ is
the order of the finite group $\Pi(G)=\Lambda/L$, and $\tau_{\lambda}(\nu_0) = \mu_{\lambda}^{-1} (\nu_0) \in X$.

\label{mwasympt}
\end{maintheo}

Next, we consider a general weight $\nu$. We have just handled the
case where $d_{N}(\nu) \sim N \nu_0$, so now we assume that
$|d_{N}(\nu)| = o(N)$, i.e. the weight lies  in the moderate
deviations region. All of the objects in the previous result
continue to make sense in this regime, but now depend on $N$.

\begin{maintheo}
Let $\lambda \in C \cap I^{*}$ be a dominant weight, and let $\nu \in NQ(\lambda)$
be a weight of the form
\[
\nu =NQ^{*}(\lambda) +d_{N}(\nu),\quad |d_{N}(\nu)|=o(N),
\]
where $|d_{N}(\nu)|$ denotes the norm of the vector $d_{N}(\nu)$ with respect to
the fixed $W$-invariant inner product on $\mf{t}^{*}$.
Assume that $m_{N}(\lambda;\nu) \neq 0$ for every sufficiently large $N$.
Then, in the notation above, we have:
\[
m_{N}(\lambda;\nu)=
(2\pi N)^{-m/2}
\frac{|\Pi(G)|e^{N\delta_{\lambda}(\nu/N)}}{\sqrt{\det A_{\lambda}(\nu/N)}}
(1+O(N^{-1})).
\]
\label{mwasympt2}
\end{maintheo}

Note that, in Theorem \ref{mwasympt2}, the point $\nu/N$ is in the interior $Q(\lambda)^{o}$ of the
polytope $Q(\lambda)$ for sufficiently large $N$, since the center of mass $Q^{*}(\lambda)$
is in the interior $Q(\lambda)^{o}$ and that the vector $d_{N}(\nu)$ is assumed to be of order $o(N)$.
Theorem \ref{mwasympt2} is regarded as an ``interpolation'' between the central limit region
and the region of moderate deviation discussed in the beginning of this section.
In fact, one can deduce Theorem \ref{repCLT} from Theorem \ref{mwasympt2}.
See Theorem \ref{lCLT} below and Proposition \ref{PRIM} in Section \ref{partition}.

\subsection{Statement of results on irreducible multiplicities.}

As we will discuss below, the multiplicities of irreducibles in
$V_{\lambda}^{\otimes N}$ can be expressed as an alternating sum
of weight multiplicities. Thus, it would be natural to expect that
one might obtain asymptotics of irreducible multiplicities from
our theorems on weight multiplicities stated above. Before stating
our result, we should mention the following result, due to Biane
\cite{B}, which gives the asymptotics of irreducible
multiplicities in the central limit region. To our knowledge, this
is the only prior result on asymptotics on pointwise
multiplicities in high tensor products.

\begin{maintheo}(Biane)  \label{BIANE}  (Th\'{e}or\`{e}me $2.2.$ in
\cite{B})
Assume that $G$ is semisimple.
For every positive integer $N$, let $NM_{\lambda}$ be the set of weights of the form
$\nu_{1}+\cdots+\nu_{N}$ with $\nu_{j} \in M_{\lambda}$.

Then, for $\mu \not \in NM_{\lambda} +\Lambda^{*}$, $a_{N}(\lambda;\mu)=0$.
For, $\mu \in NM_{\lambda}+\Lambda^{*}$ with $|\mu|
\leq C \sqrt{N}$, we have:
\[
a_{N}(\lambda;\mu) =\frac{|\Pi(G)|(\dim V_{\lambda})^{N}(\dim
V_{\mu}) \prod_{\alpha \in
\Phi_{+}}\ispa{A_{\lambda}^{-1}\alpha,\rho}} {\sqrt{\det
A_{\lambda}}(2\pi)^{m/2}N^{(\dim G)/2}} \left(
e^{-\ispa{A_{\lambda}^{-1}(\mu +\rho),\mu +\rho}/(2N)}+O(N^{-1/2})
\right),
\]
where the matrix $A$ is defined in \eqref{repA},
$m$ is the rank of $G$ and the inner
product $\ispa{\cdot,\cdot}$ is the Killing form.
\end{maintheo}

To be more precise,  in Biane \cite{B} $G$ is the  Lie group with
Lie algebra $\mf{g}$ (which is assumed to be simple in \cite{B})
such that the integral lattice of a maximal
torus is identified with the dual of the lattice
$I_{\lambda}^{*}$ generated by $M_{\lambda}$.
His quadratic form $q$ is the same as our
$A_{\lambda}$. Thus, for example, the term
$k(E)/\vol_{q}(\mf{t}/\check{Q})$ in \cite{B} is equal to our
$|\Pi(G)|/\sqrt{\det A}$ when $E=V_{\lambda}$.

 The two Theorems can be formally related by expressing the
multiplicity $a_{N}(\lambda;\mu)$  as an alternating sum of the
weight multiplicities (see Proposition \ref{Mofirred}). By
\eqref{CLTmult} one has:
\begin{gather*}
m_{N}(\lambda;\mu +\rho -w\rho) =\frac{|\Pi(G)|(\dim
V_{\lambda})^{N}}{(2\pi N)^{m/2}\sqrt{\det A_{\lambda}}}
(c_{w,N}(\lambda;\mu) +O(N^{-1/2})), \\
c_{w,N}(\lambda;\mu)= e^{-\ispa{A_{\lambda}^{-1}(\mu +\rho
-w\rho), (\mu +\rho -w\rho)}/(2N)},
\end{gather*}
Since the matrix $A_{\lambda}$ is $W$-invariant if the Lie algebra
is simple, it follows that  the quadratic form
$\ispa{A_{\lambda}^{-1}\nu,\nu}$ is a multiple of the Killing form
by some positive constant. Thus, we have
\[
\sum_{w \in W}\sgn(w) c_{w,N}(\lambda;\mu) =\frac{(\dim
V_{\mu})\prod_{\alpha \in
\Phi_{+}}\ispa{A_{\lambda}^{-1}\alpha,\mu}}{N^{d}} (1+O(N^{-1})),
\]
where $d$ is the number of the positive roots.
Therefore the alternating sum above agrees with the leading term of
Biane's formula, since $\dim G=m+2d$. However, to prove the
Theorem \ref{BIANE} in this way, one would need to prove that the
remainder similarly cancels to order $N^{-d}$ when summed over the
Weyl group, and that would be harder than the (relatively simple)
direct proof of Biane.

Although the alternating sum approach to the irreducible
multiplicities does not seem to be optimal in the central limit
region as explained above, we can deduce an asymptotic formula for
the irreducible multiplicities from Theorem \ref{mwasympt} in the
region of the strong deviations under some assumptions on the
dominant weight.

\begin{maintheo}
Let $V_{\lambda}$ be an irreducible representation of $G$ with the
highest weight $\lambda \in C$.
Let $\nu \in M_{\lambda} \cap \ol{C}$ be a dominant weight which
occurs in $V_{\lambda}$ as a weight and is assumed to lie in the
interior of the polytope $Q(\lambda)$. Then we have the following
asymptotic formula for the multiplicity $a_{N}(\lambda;N\nu)$:
\begin{equation}
a_{N}(\lambda;N\nu)=(2\pi N)^{-m/2}e^{N\delta_{\lambda}(\nu)}
\left( \frac{|\Pi(G)|\Delta (\tau_{\lambda}(\nu)/(2\pi
i))e^{-\ispa{\rho,\tau_{\lambda}(\nu)}}} {\sqrt{\det
A_{\lambda}(\nu)}} +O(N^{-1}) \right), \label{asymptF3}
\end{equation}
where $m$ is the number of simple roots, $|G_{\lambda}|$ is the
order of the finite group $G_{\lambda}=L_{\lambda}/L$. The vector
$\tau_{\lambda}(\nu) \in X$, the positive constant
$\delta_{\lambda}(\nu) >0$ and the real positive matrix
$A_{\lambda}(\nu)$ are given in \eqref{Gmomexp}, \eqref{Gconst0}
and \eqref{Gmatrix0}, respectively, and $\Delta$ is the Weyl denominator
extended to the complexification $\mf{t}^{\mb{C}}$.
\label{IMULT}
\end{maintheo}

\begin{rem}

\begin{itemize}

\item The constant $\delta_{\lambda}(\nu)$ and the matrix
$A_{\lambda}(\nu)$ are determined by the irreducible
representation $V_{\lambda}$ itself. In particular, they can be
computed by the logarithmic differential of the character of the
irreducible representation $V_{\lambda}$.

\item  The constant $\delta_{\lambda}(\nu)$ is positive under the
assumptions in Theorems \ref{mwasympt}, \ref{mwasympt2} and
\ref{IMULT}. Hence,  the multiplicities $a_{N}(\lambda;\nu)$ have
an exponential growth with respect to $N$ in the regions under
consideration.

\item It follows from  Theorem \ref{IMULT} that  the term
$\Delta (\tau_{\lambda}(\nu)/2\pi i)$ in \eqref{asymptF3} is
non-negative for such a $\nu$ as in Theorem \ref{IMULT}. We prove
this fact directly for $G=U(2)$ in Section \ref{unitary}. As the
example in Section \ref{unitary} suggests, if the dominant weight
$\nu$ is in a wall of a Weyl chamber, then the leading term in
\eqref{asymptF3} might vanish.

\end{itemize}

\end{rem}

\subsubsection{Rare events}
It should be possible to obtain further results on rare events
reminiscent of the Poisson limit law for the multinomial
distribution. Recall that
  the binomial distribution with parameter $p$  tends to a
Poisson distribution if $p \to 0$ as $N \to \infty$ with $p/N \to
C.$ Because our results allow for general coefficient weights $c $
on $S$, we believe there are  analogous results  on multiplicities
of weights  near the boundary of $Q(N \lambda).$ However, for the
sake of brevity we do not carry out the analysis of this case.

\subsubsection{Joint asymptotics}

The asymptotics of tensor products $V_{\lambda}^{\otimes N}$ as $N
\to \infty$ may be regarded as a thermodynamic limit. As recalled
in Section \ref{SYMPMOD}, the asymptotics as the highest weight
$\lambda \to \infty$ is a semiclassical limit studied by Heckman,
Guillemin-Sternberg and others.  By combining the methods of this
paper and those of Heckman et al., one could probably obtain joint
asymptotics as $N \to \infty, \lambda \to \infty$ of
multiplicities of $V_{\lambda}^{\otimes N}$. This again is
motivated by the complexity of multiplicity formulae when either
$N$ or $\lambda$ is large.

\subsection{Statement of results on lattice path multiplicities}

As mentioned above, our results on multiplicities of weights and
irreducibles are special cases of results on asymptotic counting
of lattice paths with steps in a convex lattice polytope.
Relations between lattice paths and representations have been
studied for some time, and one is  proved by Grabiner-Magyar in
(\cite{GM}). We include a proof of an adequate relation for our
purposes in Proposition \ref{Mofirred} (see also Proposition
\ref{MequalsP} for the case of weights). General and conceptually
clear relations can be derived from the path discussed in
Littelmann's expository article \cite{Lit}. We add some further
comments in Section \ref{FINAL}.

Let us now
 recall what the combinatorics of lattice paths is about:
Given a set $S \subset {\bf N}^m $ of {\it allowed steps}, an $S$-
lattice path of length $N$ from $0$ to $\beta$ is a sequence
$(v_1, \dots, v_N) \in S^N$ such that $v_1 + \cdots + v_N =
\beta$. We define the multiplicity (or partition) function of the
lattice path problem by \begin{equation} \label{LPM} {\mathcal
P}_N(\gamma) = \#\{(v_1, \dots, v_N) \in S^N: v_1 + \cdots + v_N =
\gamma\}.
\end{equation}
The set of possible endpoints of an $S$- path of length $N$ forms
a set ${\mathcal P}_{S, N},$ and we may ask how the numbers
$\mathcal{P}_{N}(\gamma)$ are distributed as $\gamma$ varies over
${\mathcal P}_{S, N}$.

It is useful (and requires no more work) to consider a somewhat
more general problem:  Let $X$ be a real vector space and let  and
$L \subset X$ be a lattice.  Also,  let $X^{*}$ and $L^{*}$ be
their duals. Let $S \subset L^{*}$ ($\# S \geq 2$) be a finite set
which satisfies the following condition:
\[
\mbox{The set}\ \ \{\beta -\beta^{\prime}\,;\,\beta,\beta^{\prime} \in S\}
\ \ \mbox{spans}\ \ X^{*}.
\]
Let $P$ be the convex hull of the finite set $S$. Let $L(S)^{*}$
be the lattice in $X^{*}$ spanned by $\{\beta
-\beta^{\prime}\,;\,\beta,\beta^{\prime} \in S\}$ over $\mb{Z}$,
and let $L(S)$ be its dual lattice. By the above assumptions, we
have $L \subset L(S)$, and the quotient $\Pi(S):=L(S)/L$ is a finite
group. For a strictly positive function $c$ on $S$, we define the
{\it weighted multiplicity of lattice paths} $\mathcal{P}_{N}^{c}$
{\it of length} $N$ {\it with weight} $c$ and {\it the set of the
allowed steps} $S$ by

\begin{equation}
\mathcal{P}_{N}^{c}(\gamma)
=\sum_{\beta_{1},\ldots,\beta_{N} \in S\,;\,\gamma =\beta_{1}+\cdots +\beta_{N}}
c(\beta_{1}) \cdots c(\beta_{N}),\quad
\gamma \in (NP) \cap L^{*}.
\label{lpcf}
\end{equation}
If  $c \equiv 1$, then $\mathcal{P}_{N}^{c}(\gamma) = {\mathcal
P}_N(\gamma)$.   In the case where $S = p \Sigma \cap \N^m$ is the
set of lattice points in the simplex of degree $p$, and the weight
function $c$ is given by $c(\beta)=\frac{p!}{\beta!(p-|\beta|)!}$,
the coefficients ${\mathcal P}^{c}_N(\gamma)$ are precisely
multinomial coefficients, and in general one may regard ${\mathcal
P}_{N}^{c}$ as a generalized multinomial coefficient. In
Proposition \ref{MequalsP}, we prove  that weight multiplicities
can be equated with weighted multiplicities of certain lattice
paths, specifically
\begin{equation} \label{MULTM} m_N(\lambda; \mu) = \pcal_N^{c_{\lambda}}(\mu - N
\lambda), \end{equation} where $\pcal_N^{c_{\lambda}}$ is a certain
weighted lattice path partition function.In Proposition
\ref{Mofirred}, we further prove that
\begin{equation} \label{MULTA} a_N(\lambda; \mu) = \Sigma_{w \in W}
\mbox{sgn}(w) \pcal_N^{c_{\lambda}}(\mu - N \lambda + \rho - w \rho).
\end{equation}

We now state our results on multiplicities of lattice paths,
following the same outline as for weight multiplicities. As in the
case of group representations, the  simplest question to consider
is the weak limit of the measure
\begin{equation} \label{COM} d \mu_{S, N}:= \frac{1}{(\# S)^N} \;
\sum_{\beta \in {\mathcal P}_{S, N}} \; {\mathcal P}_N (\beta)
\delta_{\frac{\beta}{N}}.
\end{equation}
It is well-known and easy to prove ( see Proposition \ref{lattMl})
that
\begin{equation} d \mu_{S, N} \to
\delta_{m^{*}_{S}},\;\;\;\mbox{where}\;\; m^{*}_{S}=\frac{1}{\#
S}\sum_{\beta \in S} \beta
\end{equation}
is the center of mass of the set $S$.  In the more general case of
weighted lattice paths, the center of mass $m^{*}_{S} \in P^{o}$
is given by \begin{equation} m^{*}_{S}=\frac{1}{V(S)}\sum_{\beta
\in S}c(\beta)\beta,\quad V(S)=\sum_{\beta \in S}c(\beta).
\label{centerom}
\end{equation}
We then consider the asymptotic distribution of multiplicities of
lattice paths in regions around the center point.

These refined results  involve the `moments maps',
\begin{equation} \label{MAPM}
\mu_{P}:X \to P^{o},\quad \mu_{P}(\tau)=\sum_{\beta \in S}
\frac{c(\beta)e^{\ispa{\beta,\tau}}}{\sum_{\beta' \in S}
c(\beta')e^{\ispa{\beta,\tau}}}\beta.
\end{equation}  For $x \in
P^{o}$, the interior of the polytope $P$, we define the function
$\delta_{c}(S,x)$
\begin{equation}
 \delta_{c}(S,x)=\log \left( \sum_{\beta \in
S}c(\beta)e^{\ispa{\beta -x,\tau_{P}(x)}} \right), \label{const0}
\end{equation}
and   the positive definite linear map $A_{c}(S,x):X \to X^{*}$ by
\begin{equation}  A_{c}(S,x)= \sum_{\beta \in S} \left(
\frac{c(\beta)e^{\ispa{\beta,\tau_{P}(x)}}} {\sum_{\beta' \in
S}c(\beta')e^{\ispa{\beta',\tau_{P}(x)}}} \right) \beta \otimes
\beta -x \otimes x, \quad A=A_{c}(S,m^{*}_{S}) \label{matrix0}
\end{equation}
where the diffeomorphism $\tau_{P}:P^{o} \to X$ is the inverse of
the `moment map' $\mu_{P}$.

\begin{rem}

It should be noted that the constant $\delta_{c}(S,\alpha)$
defined in \eqref{const0} depends on the choice of the weight
function $c$. In fact this constant can be negative if we choose
the weight function $c$ small enough. However, if $c$ takes
positive integer values, then it turns out that the constant
$\delta_{c}(S,\alpha)$ is positive. See {\it Remark} after the
proof of Theorem \ref{lpasympt} in Section \ref{partition}.

\end{rem}

\subsubsection{Central limit region}

Our first result on lattice paths concerns the central limit region where $\gamma
=Nm^{*}_{S} + d_{N}(\gamma),\quad d_{N}(\gamma)=O(N^{s})$ for a
variety of $s < 1.$

\begin{maintheo}
Let $0 \leq s < 1$.
Let $\gamma$ be a lattice point such that $\mathcal{P}_{N}^{c}(\gamma) \neq 0$ for every
sufficiently large $N$, and assume also that $\gamma$ has the form
\begin{equation}
\gamma =Nm^{*}_{S} + d_{N}(\gamma),\quad
d_{N}(\gamma)=O(N^{s}).
\label{order}
\end{equation}
Then we have
\begin{equation}
\mathcal{P}_{N}^{c}(\gamma)=
(2\pi N)^{-m/2}\frac{|\Pi(S)|e^{N\delta_{c}(S,\gamma/N)}}{\sqrt{\det A}}(1+O(N^{-(1-s)})).
\label{genasympt}
\end{equation}
Furthermore, if $0 \leq s \leq 2/3$ and $d_{N}(\gamma)=o(N^{s})$, we have
\begin{equation}
\mathcal{P}_{N}^{c}(\gamma)
=(2\pi N)^{-m/2}\frac{|\Pi(S)|V(S)^{N}e^{-\ispa{A^{-1}d_{N}(\gamma),d_{N}(\gamma)}/(2N)}}
{\sqrt{\det A}}(1+\varepsilon_{N}),
\label{localclt}
\end{equation}
where
\[
\varepsilon_{N}=
\left\{
\begin{array}{ll}
O(N^{-(1-s)}) & \mbox{for} \quad 0 \leq s \leq 1/2, \\
o(N^{3s-2}) & \mbox{for} \quad 1/2 < s \leq 2/3.
\end{array}
\right.
\]
\label{lCLT}
\end{maintheo}

\subsubsection{Large deviations region}

We now assume that $d_N $ is of order $N$.

\begin{maintheo}
Let $\alpha$ be a lattice point in $S$ which is assumed to lie in
the interior of the polytope $P$. Then, for every $f \in
L(S)^{*}$, we have
\begin{equation}
\mathcal{P}_{N}^{c}(N\alpha +f)= (2\pi
N)^{-m/2}\frac{|\Pi(S)|e^{-\ispa{f,\tau_{P}(\alpha)}+N\delta_{c}(S,\alpha)}}{\sqrt{\det
A_{c}(S,\alpha)}} (1+O(N^{-1})), \label{asymptF1}
\end{equation}
where $|\Pi(S)|$ denotes the order of the finite group
$\Pi(S)=L(S)/L$. The exponent $\delta_{c}(S,\alpha)$ is positive if
$c(\alpha) \geq 1$. \label{lpasympt}
\end{maintheo}

Our analysis starts from the fact that
\[
\mathcal{P}_{N}(\gamma) = \chi_S(u)^N |_{u^{\gamma}},
\]
where $\chi_S(u)^N |_{u^{\gamma}}$ denotes the coefficient of the
monomial  $u^{\gamma}$ in the $N$-th power of the {\it admissible
step character},
\begin{equation} \chi_S(u) = \sum_{\alpha \in S} u^{\alpha}
\end{equation}
We apply a steepest descent argument  to an integral
representation of $\mathcal{P}_{N}^{c}(\gamma)$ (see
\eqref{integral1} in Section \ref{partition}). Our basic reference
for the stationary phase for complex phase functions is
(\cite{Ho}).

\subsection{Organization} We first prove the results on lattice
paths, Theorems \ref{lCLT} and \ref{lpasympt}, in Section
\ref{partition}. We then deduce the main  results on
multiplicities, Theorems \ref{repCLT}, \ref{BIANE},
\ref{mwasympt}, \ref{mwasympt2}, \ref{IMULT}, in Section
\ref{GREP}. In that section, we also review the relation between
multiplicities of weights and lattice paths. In Section
\ref{unitary}, we illustrate the results for some representations
of $U(m)$ with $m=2$. In Section \ref{FINAL}, we make some final
comments on the connections between lattice paths and weight
multiplicites and on the symplectic model for tensor product
multiplicities.
\medskip

\noindent{{\bf Acknowledgments}} This paper originated as a
by-product  of our joint work \cite{TSZ1} with B. Shiffman on the
\szego kernel of a toric variety and its applications, where the
partition function of a certain lattice path problem appeared. The
probabilistic background became clearer in discussions on the
local central limit theorem and large deviations with  Y. Peres.

This paper was written during a stay of the first author at the
Johns Hopkins University on a JSPS fellowship, and during the
visit of both authors at MSRI at the program on semiclassical
analysis.

\section{Asymptotics of the number of Lattice paths}
\label{partition}

Let $X$ be a finite dimensional real vector space of dimension $m$, and let
$L$ be a lattice in $X$. Let $X^{*}$ and $L^{*}$ be, respectively,
the dual vector space of $X$ and the dual lattice of $L$.
Let $S \subset L^{*}$ be a finite set such that $\# S \geq 2$, and set
\begin{equation}
D(S):=\{\beta-\beta^{\prime} \in L^{*}\,;\,\beta,\beta^{\prime} \in S\}.
\label{diffset}
\end{equation}
We assume that
\begin{equation}
\lspan_{\mb{R}}D(S)=X^{*}.
\label{assumpt1}
\end{equation}
Let $P=P_{S}$ be the convex hull of $S$, which is an integral polytope in $X^{*}$.
Let
\[
c:S \to \mb{R}_{>0}
\]
be a strictly positive function on $S$.
Our aim in this section is to investigate the asymptotics of
the number of the lattice paths $\mathcal{P}_{N}^{c}(\gamma)$
for the lattice point $\gamma$ in various regions
(central limit region, regions of moderate and strong deviations
discussed in Introduction) as $N \to \infty$.

We introduce the {\it weighted character} (or the {\it weighted $S$-character})
with the weight function $c$ defined by
\begin{equation}
k(w):=\sum_{\beta \in S}c(\beta)e^{\ispa{\beta, w}},
\quad w \in X^{\mb{C}}:=X \otimes \mb{C},
\label{pchar1}
\end{equation}
which is considered as a function on $X^{\mb{C}} =X \otimes \mb{C}$.
Here, and in what follows, a functional $f \in X^{*}$ is considered as a $\mb{C}$-linear
functional on $X^{\mb{C}}$. We fix a primitive basis for the lattice $L$, which is also
considered as a fixed basis for $X$.
Note that, for $\tau \in X$,
the function $\varphi \mapsto k(\tau +i\varphi)$ is a smooth
function on the torus ${\bf T}^{m}:=X/(2\pi L)$, since we have assumed $S \subset L^{*}$.
The fixed basis in $L$ defines a Lebesgue measure on $X$, and hence on $X^{\mb{C}}$,
normalized so that $\vol ({\bf T}^{m})=(2\pi)^{m}$.
We also fix an inner product on $X$ which has the fixed basis for $L$
as an orthonormal basis, and we denote by $|\varphi|$ the norm of $\varphi \in X$
with respect to this inner product.

It is clear that the $N$-th power of the function $k(w)$ is given by
\begin{equation}
k(w)^{N}=\sum_{\gamma \in (NP) \cap L^{*}}
\mathcal{P}_{N}^{c}(\gamma)e^{\ispa{\gamma,w}}.
\label{Nthpower}
\end{equation}
Therefore, the lattice paths counting function $\mathcal{P}_{N}^{c}$ has the following
integral expression:
\begin{equation}
\mathcal{P}_{N}^{c}(\gamma)
=\frac{1}{(2\pi)^{m}}\int_{{\bf T}^{m}}
e^{-i\ispa{\gamma,\varphi}}k(i\varphi)^{N}\,d\varphi.
\label{integral1}
\end{equation}

To begin with, we shall consider the simplest case, that is,
consider the problem how the numbers of lattice paths
with endpoints varying in $NP \cap L^{*}$ distributes. This would be expressed as
the weak limit of the measure defined by the following:
\begin{equation}
m_{S,N}:=\frac{1}{V(S)^{N}}\sum_{\gamma \in NP \cap L^{*}}
\mathcal{P}_{N}^{c}(\gamma)\delta_{\gamma/N},\quad
V(S):=k(0)=\sum_{\beta \in S}c(\beta).
\label{measurelatt}
\end{equation}
Noting that $\mathcal{P}_{1}^{c}(\gamma)=c(\gamma)$ $(\gamma \in S)$, we have
\[
V(S)^{N} =\sum_{\gamma \in NP \cap L^{*}}\mathcal{P}_{N}^{c}(\gamma),
\]
which shows that the measure $m_{S,N}$ is a probability measure.
The following proposition will be used to prove Theorem \ref{SIMPLE} in the next section.

\begin{prop}
The probability measure $m_{S,N}$ tends weakly to the Dirac measure $\delta_{m^{*}_{S}}$ at
the point $m^{*}_{S} \in P$ given in \eqref{centerom}.
\label{lattMl}
\end{prop}

\begin{proof}
It suffices to show that the Fourier transform (characteristic function) $\widehat{m_{S,N}}(\varphi)$
of the probability measure $m_{S,N}$ converges to the Fourier transform of the Dirac measure
$\delta_{m^{*}_{S}}$ at the point $m^{*}_{S}$ for every $\varphi \in X$.
The Fourier transform of $\delta_{m^{*}_{S}}$ is given by $\varphi \mapsto e^{-i\ispa{m^{*}_{S},\varphi}}$.
By \eqref{integral1}, the Fourier transform of $m_{S,N}$ is given by
\[
\widehat{m_{S,N}}(\varphi)
=\left[
\frac{k(-i\varphi/N)}{V(S)}
\right]^{N},\quad \varphi \in X.
\]
Thus we need to show that $\widehat{m_{S,N}}(\varphi) \to e^{-i\ispa{m^{*}_{S},\varphi}}$ as $N \to \infty$.
Since $\widehat{m_{S,N}}(0)=1$, we can choose a compact neighborhood $U$ of the origin in $X$ such that
a branch of the logarithm $\log \widehat{m_{S,N}}(\varphi)$ exists for $\varphi \in U$.
For any $\varphi \in X$ we take $N$ large enough so that $\varphi/N \in U$. Then, a Taylor expansion
at the origin gives
\[
e^{N\log\widehat{m_{S,N}}(\varphi/N)}
=e^{-i\ispa{m^{*}_{S},\varphi}+N^{-1}R_{N}(\varphi)},
\]
where $R_{N}(\varphi)$ is bounded on compact sets uniformly in $N$.
Therefore, we have
$\widehat{m_{S,N}}(\varphi) \to e^{-i\ispa{m^{*}_{S},\varphi}}$ as $N \to \infty$.
\end{proof}

Next, although our main purpose of this section is to find various aspects
of asymptotics of the weighted number (multiplicity) of $S$-lattice paths,
we shall give some account on the large deviation principle for the sequence of
probability measures.

\begin{prop}
The sequence of measures $\{m_{S,N}\}$ satisfies the large deviation principle
with the rate function given by
\begin{equation}
I_{S}(x)=\sup_{\tau \in X}
\left \{
\ispa{\tau,x}-\log
(k(\tau)/V(S))
\right \}.
\label{latRATE}
\end{equation}
\label{latLDP}
\end{prop}

\begin{proof}
We apply Cram\'{e}r's theorem (\cite{DZ}, Theorem 2.2.30). We
shall recall the setting-up for the Cram\'{e}r's theorem. Let
$X_{j}$ ($j=1,2,\ldots$) be a sequence of independent identically
distributed $m$-dimensional random vectors on a probability space
with $X_{1}$ distributed according to the probability measure
$\mu$ on $\mb{R}^{m}$. Let $m_{N}$ be the distribution
(probability measure) for the empirical means
$S_{N}:=\frac{1}{N}\sum_{j=1}^{N}X_{j}$. Then, Cramer's theorem
states that the sequence of measures $\{m_{N}\}$ satisfies the LDP
with the rate function
\[
I(x)=\sup_{\tau \in \mb{R}^{m}}
\{\ispa{\tau,x}-\Lambda(\tau)\},\quad
\Lambda(\tau)={\bf E}(e^{\ispa{\tau,X_{1}}})
\]
if $\Lambda(\tau) < \infty$ for every $\tau \in \mb{R}^{m}$.
In our case, We take the probability space
$\mathcal{E}:=P \times \cdots$ (infinite product of the polytope $P$),
and the probability measure $m_{S} \times \cdots$ on $\mathcal{E}$.
The random variable $X_{j}$ is the projection onto the $j$-th factor.
Then, it is easy to see that $\Lambda(\tau)=\log (k(\tau)/V(S))$, and the
push-forward of the measure $m_{S} \times \cdots$ by the empirical mean
$S_{N}=\frac{1}{N}\sum_{j=1}^{N}X_{j}$ is nothing but $m_{S,N}$.
Therefore, the assertion is a direct consequence of Cram\'{e}r's theorem stated above.
\end{proof}

Proposition \ref{lattMl} suggests that the number of lattice paths would have a `peak' at
the center of mass (although, in general, the center of mass might not be in the lattice $L^{*}$).
Thus, it is natural to ask that how the lattice paths counting function $\mathcal{P}_{N}^{c}(\gamma)$
behave with the distance between $\gamma$ and the center of mass getting large.
But, when $N$ becomes large, the possible end points of the $S$-lattice paths is in the polytope $NP$,
and the center of mass of $NP$ is $Nm^{*}_{S}$ where $m^{*}_{S}$ is the center of mass of $P$ defined in \eqref{centerom}.
Thus it is natural to consider the behavior of $\mathcal{P}_{N}^{c}(\gamma)$ when the distance between $\gamma$ and
$Nm^{*}_{S}$ varies.

Our next aim in this section is to prove Theorems \ref{lCLT} and
\ref{lpasympt} which corresponds respectively the the case where
$\gamma$ is in the central limit region (and the region of
moderate deviations) and the region of the strong deviations.

\subsection{Proof of Theorem \ref{lpasympt} }

First we shall prove Theorem \ref{lpasympt}. To prove Theorem \ref{lpasympt}, we
need to prepare notation.

Let $\exp:X \to {\bf T}^{m}:=X/(2\pi L)$ be the exponential map, {\it i.e.},
the canonical projection.
Since the set of differences $D(S)$ defined in \eqref{diffset} spans $X^{*}$,
it spans a lattice, $L(S)^{*}$, in $X^{*}$ of maximal rank over $\mb{Z}$:
\begin{equation}
L(S)^{*}=\lspan_{\mb{Z}}D(S) \subset L^{*},
\end{equation}
and its dual lattice in $X$ is denoted by $L(S)$.
We have $L(S)^{*} \subset L^{*}$, and hence $L \subset L(S)$. Both of the lattices is of maximal rank.
Thus, the quotient group $\Pi(S):=L(S)/L$ is a finite group.
The finite group $\Pi(S)$ is the covering transformation group of the surjective homomorphism
\begin{equation}
\pi_{S}:{\bf T}^{m} \to T(S):=X/(2\pi L(S)),\quad
\pi_{S}(\exp \varphi)=\exp_{S}(\varphi),
\label{exacts1}
\end{equation}
where $\exp_{S}:X \to T(S)$ denotes the canonical projection.

\begin{rem}
If we begin with a polytope $P$, the function $c$ above should be a non-negative function
on $P \cap L^{*}$. In this case, the corresponding finite set $S$ should be the support of the function $c$.
Thus, the support $S$ of the function $c$ is assumed to satisfy \eqref{assumpt1}.
If the set $D(S)$ defined in \eqref{diffset} spans the lattice $L^{*}$ over $\mb{Z}$, then the corresponding
torus $T(S)$ coincides with the original torus ${\bf T}^{m}$, and hence $\Pi(S)=\{1\}$.
\end{rem}

\begin{lem}
For any fixed vector $\tau \in X$, we denote $k_{\tau}(\exp \varphi):=k(\tau +i\varphi)$, which
is considered as a function on ${\bf T}^{m}$, where the function $k$ on $X^{\mb{C}}$ is given in \eqref{pchar1}.
Then we have $|k_{\tau}(\exp \varphi)| \leq k(\tau)$. The equality holds exactly on the kernel of the
homomorphism $\pi_{S}:{\bf T}^{m} \to T(S)$:
\[
\{t \in {\bf T}^{m}\,;\,|k_{\tau}(t)| =k(\tau)\}=\ker \pi_{S} \cong \Pi(S).
\]
In particular, the set in the left hand side is finite.
\label{equalone}
\end{lem}

\begin{proof}
The inequality $|k_{\tau}(\exp \varphi)| \leq k(\tau)$ follows from the Cauchy-Schwarz inequality.
It is easy to see that the condition $|k_{\tau}(\exp \varphi)|=k(\tau)$ on $\varphi \in X$ is equivalent to
the following:
\[
\ispa{\beta-\beta^{\prime},\varphi} \in 2\pi \mb{Z},\quad \beta,\beta^{\prime} \in S.
\]
Since $L(S)^{*}=\lspan_{\mb{Z}}D(S)$, this condition is equivalent to say that $\varphi \in 2\pi L(S)$.
This completes the proof.
\end{proof}

Note that the function $k(w)=k(\tau +i\varphi)$ is holomorphic in $w=\tau +i \varphi \in X^{\mb{C}}$,
and is $2\pi L$-periodic with respect to the variable $\varphi \in X$.
Therefore, we can deform the contour of the integral in \eqref{integral1},
and hence, by setting $\gamma =N\alpha +f$ in \eqref{integral1}, we can write
\begin{equation}
\mathcal{P}_{N}^{c}(N\alpha +f)
=\frac{e^{-\ispa{f,\varphi}}}{(2\pi)^{m}}[k(\tau)e^{-\ispa{\alpha,\tau}}]^{N}
\int_{{\bf T}^{m}}e^{-iN\ispa{\alpha,\varphi}}
\left [
\frac{k(\tau +i\varphi)}{k(\tau)}
\right]^{N}e^{-i\ispa{f,\varphi}}\,d\varphi,
\label{integral2}
\end{equation}
where $\tau \in X$ is arbitrary. (Note that $k(\tau) >0$ for $\tau \in X$.)
To choose a suitable $\tau \in X$, we need to find the point where the function
$k(\tau)e^{-\ispa{\alpha,\tau}}$ attains its minimum.
To describe the critical points of this function, we define a map $\mu_{P}: X \cong \mb{R}^{m} \to P^{o}$ by
\begin{equation}
\mu_{P}(\tau):=\partial_{\tau}\log k(\tau)
=\frac{1}{\sum_{\beta \in S}c(\beta)e^{\ispa{\beta,\tau}}}
\sum_{\beta \in S}c(\beta)e^{\ispa{\beta,\tau}}\beta.
\label{moment}
\end{equation}
The map $\mu_{P}$ defined above is an analogue of the moment map for
a Hamiltonian torus action on toric manifolds.
Thus we call the map $\mu_{P}$ the {\it moment map}.
Since the set $D(S)$ of differences of vectors in the finite set $S$
spans the whole space $X^{*}$ (over $\mb{R}$), the elements in $S$ are not contained simultaneously in
any affine hyperplane in $X^{*}$. It is well-known (p.~83 in \cite{Fu}) that the moment map $\mu_{P}$
defines a (real analytic) diffeomorphism between the vector space $X$ and the interior $P^{o}$
of the polytope $P$.

We denote the inverse of the moment map $\mu_{P}$ by $\tau_{P}=\tau_{P}(x):P^{o} \to X$.
Then, for every $\alpha \in P^{o}$, we have $\mu_{P}(\tau_{P}(\alpha))=\alpha \in P^{o}$.

We note that the center of mass $m^{*}_{S}$ is the value of the moment map at the origin:
$\mu_{P}(0)=m^{*}_{S}$, $\tau_{P}(m^{*}_{S})=0$.
The differential of the moment map $\mu_{P}:X \to P^{o}$ defines the following linear transform
$A(\tau):X \to X^{*}$.
\[
A(\tau):=\sum_{\beta \in S}\frac{c(\beta)e^{\ispa{\tau,\beta}}}{k(\tau)}\beta \otimes \beta
-\mu_{P}(\tau) \otimes \mu_{P}(\tau),\quad \tau \in X, \quad
A:=A(0).
\]

\begin{lem}
We set
\begin{equation}
f_{\alpha}(\tau):=\log k(\tau)-\ispa{\alpha,\,\tau},\quad
\tau \in X,
\label{aux11}
\end{equation}
so that $e^{f_{\alpha}(\tau)}=k(\tau)e^{-\ispa{\alpha,\,\tau}}$.
Then the Hessian of the function $f_{\alpha}$, which is given by $A(\tau)$,
is a positive definite for every $\tau \in X$.
The vector $\tau_{P}(\alpha)$ is the unique critical point
of the function $f_{\alpha}$. In fact, we have
\[
f_{\alpha}(\tau) \geq f_{\alpha}(\tau_{P}(\alpha)),\quad \tau \in X
\]
with equality holds only at $\tau=\tau_{P}(\alpha)$.
\label{critical1}
\end{lem}

\begin{proof}
It is straight forward to show that
\begin{equation}
\partial f(\tau)=\mu_{P}(\tau) -\alpha,\quad
A(\tau)=\partial^{2}f(\tau).
\label{firstderiv}
\end{equation}
Although one can prove the positivity of the map $A(\tau)$ for every $\tau \in X$
by exactly the same argument as in \cite{SZ2}, we give a proof of it for completeness.
For each $\beta \in S$, we set $m_{\beta}(\tau):=c(\beta)e^{\ispa{\beta,\tau}}/k(\tau)$ so that
$\sum_{\beta \in S}m_{\beta}(\tau)=1$.
We define a probability measure $\nu_{S}^{\tau}$ on $X^{*}$ supported on $S$, depending on $\tau \in X$, by
$d\nu_{S}^{\tau}=\sum_{\beta \in S}m_{\beta}(\tau)\delta_{\beta}$,
where $\delta_{\beta}$ denotes the Dirac measure at $\beta$.
Then, for any vector $x \in X$, we have
\[
\ispa{A(\tau)x,x}=
\int_{X^{*}}g_{x}(v)^{2}\,d\nu_{S}^{\tau}(v)-
\left|
\int_{X^{*}}g_{x}(v)\,d\nu_{S}^{\tau}(v)
\right|^{2} \geq 0,
\]
where $g_{x}$ is a linear function on $X^{*}$ defined by $g_{x}(v)=\ispa{v,x}$,
$v \in X^{*}$, $x \in X \cong \mb{R}^{m}$.
The equality in the above holds if and only if $g_{x}$ is constant on $S$.
In such a case, the function $g_{x}$ is zero on $D(S)$, since $g_{x}$ is linear.
Thus, by the assumption \eqref{assumpt1}, $g_{x}$ is zero on $X^{*}$, and which implies $x=0$.
This shows that $A(\tau)$ is positive definite for any $\tau \in X$.

By \eqref{firstderiv}, the vector $\tau_{P}(\alpha)$ is the unique critical point of the function $f_{\alpha}$,
since the map $\mu_{P}:X \to P^{o}$ is a diffeomorphism.
A Taylor expansion at $\tau=\tau_{P}(\alpha)$ for the function $f_{\alpha}$ gives
\[
f_{\alpha}(\tau)=f_{\alpha}(\tau_{P}(\alpha))+
\int_{0}^{1}(1-t)
\ispa{A(\tau_{P}(\alpha)+t(\tau-\tau_{P}(\alpha)))(\tau-\tau_{P}(\alpha)),\,\tau-\tau_{P}(\alpha)}\,dt.
\]
Since $A(\tau)$ is positive definite,
the last integral is non-negative, and equals zero if and only if $\tau=\tau_{P}(\alpha)$.
This completes the proof.
\end{proof}

It should be noted that the constant $\delta_{c}(S,\alpha)$ and the matrix $A_{c}(S,\alpha)$
defined by \eqref{const0}, \eqref{matrix0}
in Theorem \ref{lpasympt} can be written as
\begin{gather}
A_{c}(S,\alpha)=A(\tau_{P}(\alpha)),
\label{Amatrix}\\
\delta_{c}(S,\alpha)=f_{\alpha}(\tau_{P}(\alpha)).
\label{Const}
\end{gather}
Hence the matrix $A_{c}(S,\alpha)$ is real symmetric and positive definite.
It should be noted that the function $\delta_{c}(S,x)$ on $P^{o}$ defined in \eqref{const0} satisfies
\begin{equation}
\delta_{c}(S,x)=\log (V(S))-I_{S}(x), \quad x \in P^{o},
\label{DELTAvsRATE}
\end{equation}
where the function $I_{S}$ is the rate function defined by \eqref{latRATE}.

We choose the vector $\tau \in X$ in \eqref{integral2} as $\tau=\tau_{P}(\alpha)$.
Recall that, by Lemma \ref{equalone}, the absolute value of the integrand in \eqref{integral2} equals one
precisely on the set $\ker \pi_{S} \subset {\bf T}^{m}$, where $\pi_{S}:{\bf T}^{m} \to T(S)$ is a homomorphism.
The set $\ker \pi_{S}$ is a subgroup in ${\bf T}^{m}$ and isomorphic to $\Pi(S)=L(S)/L$, which is a finite group.
For each $g \in \ker \pi_{S} \cong \Pi(S)$, we take a representative $\varphi_{g} \in X$
so that $g=\exp \varphi_{g}$.
Let $V_{g} \subset U_{g}$ be open neighborhoods of the vector $\varphi_{g} \in X$
such that $U_{g} \cap \ker \pi_{S} =\{g\}$ and $\ol{V_{g}} \subset U_{g}$, and a branch of the logarithm
\[
\log \left(
\frac{k(\tau_{P}(\alpha)+i\varphi)}{k(\tau_{P}(\alpha))}
\right)
\]
exists on each of $U_{g}$. We choose a constant $c>0$ so that
\[
|k(\tau_{P}(\alpha)+i\varphi)/k(\tau_{P}(\alpha))| \leq e^{-c} \ \
\mbox{for}\ \
\exp \varphi \in {\bf T}^{m} \setminus \bigcup_{g \in \ker \pi_{S}}V_{g}.
\]
Let $\chi_{g}$ be a smooth function on $X$ supported in the open set $U_{g}$ and equals one near $V_{g}$.
Then we can write the integral \eqref{integral2} in the following form:
\begin{equation}
\mathcal{P}_{N}^{c}(N\alpha +f)
=\frac{e^{N\delta_{c}(S,\alpha) -\ispa{f,\tau_{P}(\alpha)}}}{(2\pi)^{m}}
\left (
\sum_{g \in \ker \pi_{S}}
\int_{{\bf T}^{m}}e^{N\Phi_{\alpha,g}(\varphi)}\chi_{g}(\varphi)e^{-i\ispa{f,\varphi}}\,d\varphi
+O(e^{-Nc})
\right),
\label{integral3}
\end{equation}
where the phase function $\Phi_{\alpha,g}(\varphi)$ is given by
\[
\Phi_{\alpha,g}(\varphi)=
\log \left(
\frac{k(\tau_{P}(\alpha)+i\varphi)}{k(\tau_{P}(\alpha))}
\right)-i\ispa{\alpha,\,\varphi}.
\]

By definition, the vectors $\varphi_{g}$ are in $2\pi L(S)$.
This implies that $\ispa{\beta-\beta^{\prime},\varphi_{g}}$ is $2\pi$ times an integer
for any $\beta,\beta^{\prime} \in S$. Therefore, the complex number
\begin{equation}
h(g):=e^{i\ispa{\beta,\varphi_{g}}} \in U(1),\quad \beta \in S,\quad g \in \ker \pi_{S}
\label{holon}
\end{equation}
does not depend on the choice of $\beta \in S$ and $\varphi_{g} \in \exp^{-1}(g) \subset X$.
Furthermore, we have
\begin{equation}
k(\tau+i\varphi_{g})=h(g)k(\tau),\quad
(\partial_{\varphi} k)(\tau +i\varphi_{g})=ih(g)(\partial k)(\tau),\quad
\tau \in X.
\label{equivariant1}
\end{equation}

\begin{lem}
For each $g \in \ker \pi_{S} \cong \Pi(S)$, we set
\[
C_{g}:=\{\varphi \in U_{g}\,;\,\Re\Phi_{\alpha,g}(\varphi)=0,\
\partial_{\varphi}\Phi_{\alpha,g}(\varphi)=0\}.
\]
Then we have $C_{g}=\{\varphi_{g}\}$.
Furthermore, we have
\[
e^{N\Phi_{\alpha,g}(\varphi_{g})}=h(g)^{N}e^{-iN\ispa{\alpha,\varphi_{g}}}, \quad
{\rm Hess}(\Phi_{\alpha,g})(\varphi_{g})=-A_{c}(S,\alpha).
\]
\label{critical2}
\end{lem}

\begin{proof}
That the real part of the phase function $\Phi_{\alpha,g}$ is less than or equal to zero follows from
the Cauchy-Schwarz inequality, since we have the obvious identity
\[
\Re\Phi_{\alpha,g}(\varphi)=
\log \left(
\frac{|k(\tau_{P}(\alpha)+i\varphi)|}{k(\tau_{P}(\alpha))}
\right).
\]
By the above identity and Lemma \ref{equalone}, $\Re\Phi_{\alpha,g}(\varphi)=0$ for $\varphi \in U_{g}$
if and only $\varphi=\varphi_{g}$. Thus the critical set $C_{g}$ is empty or
consists of the point $\varphi_{g}$. By \eqref{equivariant1}, we have
\[
(\partial_{\varphi}\Phi_{\alpha,g})(\varphi_{g})=
i\left[
\frac{(\partial k)(\tau_{P}(\alpha)+i\varphi_{g})}{k(\tau_{P}(\alpha))}
-\alpha
\right]
=i[\mu_{P}(\tau_{P}(\alpha))-\alpha]=0,
\]
which shows $C_{g}=\{\varphi_{g}\}$.
The rest of the assertion can be proved by a similar calculation by using the identity \eqref{equivariant1}.
\end{proof}

\vspace{5pt}

\noindent{{\it Completion of proof of Theorem \ref{lpasympt}.}}\hspace{3pt}
Let $\alpha \in S$ and $f \in L(S)^{*}$.
We set
\[
I_{g}:=\int_{{\bf T}^{m}}e^{N\Phi_{\alpha,g}(\varphi)}\chi_{g}(\varphi)
e^{-i\ispa{f,\varphi}}\,d\varphi.
\]
so that, by \eqref{integral3}, the lattice paths counting function
$\mathcal{P}_{N}^{c}(N\alpha +f)$ is written as
\[
\mathcal{P}_{N}^{c}(N\alpha +f)
=\frac{e^{N\delta_{c}(S,\alpha) -\ispa{f,\tau_{P}(\alpha)}}}{(2\pi)^{m}}
\left(
\sum_{g \in \ker \pi_{S}}I_{g}
+O(e^{-cN})
\right)
\]
for some constant $c>0$.
To obtain an asymptotic estimate for the integral $I_{g}$,
we shall use the method of stationary phase with a complex phase function.
In fact, by Lemma \ref{critical2} and Theorem 7.7.5 in \cite{Ho}, we have
\begin{equation}
I_{g}=
\left(
\frac{N}{2\pi}
\right)^{-m/2}
\frac{e^{N\Phi_{\alpha,g}(\varphi_{g}) -i\ispa{f,\varphi_{g}}}}{\sqrt{\det A_{c}(S,\alpha)}}
(1+O(N^{-1})).
\label{stationary}
\end{equation}
Since $f \in L(S)^{*}$ and $\varphi_{g} \in 2\pi L(S)$, $\ispa{f,\varphi_{g}}$ is $2\pi$ times
an integer. Furthermore, we have assumed that $\alpha \in S$.
Therefore, by Lemma \ref{critical2} and the definition of $h(g) \in U(1)$, we have
\[
e^{N\Phi_{\alpha,g}(\varphi_{g}) -i\ispa{f,\varphi_{g}}}=
h(g)^{N}e^{-i\ispa{N\alpha+f,\varphi_{g}}}=1,
\]
which shows the asymptotic formula \ref{asymptF1}.
As for the constant $\delta_{c}(S,\alpha)$, by taking the exponential $e^{\delta_{c}(S,\alpha)}$,
it is easy to prove that $\delta_{c}(S,\alpha)>0$ if $c(\alpha) \geq 1$.
\hfill\qedsymbol

\vspace{5pt}

\begin{rem}
The constant $\delta_{c} (S,\alpha)$ can be negative.
To be precise, we set $c=\max_{\beta \in S} c(\beta)$, and $f=0$.
Then $\mathcal{P}_{N}^{c}(N\alpha) \leq c^{N} \mathcal{P}_{N}^{1}(N\alpha)$,
where $\mathcal{P}_{N}^{1}(N\alpha)$
is the number of (non-weighted) lattice paths
\[
\mathcal{P}_{N}(N\alpha)=\sharp
\{
(\beta_{1},\ldots,\beta_{N}) \in S^{N}\,;\,
N\alpha =\beta_{1} +\cdots +\beta_{N}
\}.
\]
Thus if $c < e^{-\delta_{1}(S;\alpha)}$, then $\mathcal{P}_{N}^{c}(N\alpha)$ decays exponentially.
This proves that if $c(\beta) < e^{-\delta_{1}(S;\alpha)}$, then we have $\delta_{c}(S,\alpha) <0$.
\end{rem}

\vspace{5pt}

\subsection{Proof of Theorem  \ref{lCLT} }

Next, we shall prove Theorem \ref{lCLT}.
The same method as in the proof of Theorem \ref{lpasympt} will show the following

\begin{prop}
Let $\gamma =Nm^{*}_{S} +d_{N}(\gamma)$ be a lattice point in $L^{*}$ with $d_{N}(\gamma)=o(N)$.
Assume that $\mathcal{P}_{N}^{c}(\gamma) \neq 0$ for every sufficiently large $N$.
Then, we have
\begin{equation}
\mathcal{P}_{N}^{c}(\gamma)=
(2\pi N)^{-m/2}
\frac{|\Pi(S)|e^{N\delta_{c}(S,\gamma/N)}}{\sqrt{\det A_{c}(S,\gamma/N))}}(1+O(N^{-1})).
\label{prim}
\end{equation}
\label{PRIM}
\end{prop}

\begin{proof}
The proof is almost the same as the proof of Theorem \ref{lpasympt},
so we give its proof briefly. As in \eqref{integral3}, we can write
\[
\mathcal{P}_{N}^{c}(\gamma)=
\frac{e^{N\delta_{c}(S,\gamma/N)}}{(2\pi)^{m}}
\left(
\sum_{g \in \ker \pi_{S}}\int_{{\bf T}^{m}}
e^{N\Psi_{N,\gamma}(\varphi)}\chi_{g}(\varphi)\,d\varphi
+O(e^{-Nc})
\right)
\]
for some constant $c>0$, where
with the phase function $\Psi_{N,\gamma}$ is given by
\[
\Psi_{N,\gamma}(\varphi)=
\log \left[
\frac{k(\tau_{P}(\gamma/N)+i\varphi)}{k(\tau_{P}(\gamma/N))}
\right]
-i\ispa{\gamma/N,\varphi}.
\]
Here, it should be noted that $\delta_{c}(S;\gamma/N)=\log k(\tau_{P}(\gamma/N))-\ispa{\gamma/N,\tau_{P}(\gamma/N)}$.
The phase function $\Psi_{\gamma,N}$ satisfies $\Re \Psi_{\gamma,N} \leq 1$, and the point $\varphi_{g}$
is the only critical point on the support of $\chi_{g}$.
The Hessian of $\Psi_{\gamma,N}$ at $\varphi_{g}$ is $-A(\tau_{P}(\gamma/N))=-A_{c}(S,\gamma/N)$.
Although the phase $\Psi_{\gamma,N}$ depends on $N$, it is directly shown that its $C^{4}$-norm
on the support of the cut-off function $\chi_{g}$ is bounded in $N$.
Since $s <1$ and $\tau_{P}(m^{*}_{S})=0$,
we have $\gamma/N \to m^{*}_{S}$ as $N \to \infty$
and hence $A(\tau_{P}(\gamma/N)) \to A(\tau_{P}(m^{*}_{S}))=A$ as $N \to \infty$.
This shows that the norm of $A(\tau_{P}(\gamma/N))$ is bounded from below uniformly in $N$.
We have assumed that $\mathcal{P}_{N}^{c}(\gamma) \neq 0$ for every sufficiently large $N$, and hence
we have $h(g)^{N}e^{-i\ispa{\gamma,\varphi_{g}}}=1$ for any $g \in \ker \pi_{S}$ for every sufficiently large $N$.
This shows that $e^{N\Psi_{\gamma,N}(\varphi_{g})}=1$. Therefore, the assertion follows from Theorem 7.7.5 in \cite{Ho}.
\end{proof}

\noindent{{\it Completion of proof of Theorem \ref{lCLT}.}}\hspace{3pt}
First, note that we have set $A=A(0)$.
Thus, we have $\sqrt{\det A(\tau)}=\sqrt{\det A}(1+O(|\tau|))$ near $\tau =0$.
Noting $\gamma/N -m^{*}_{S} =N^{-1}d_{N}(\gamma)=O(N^{-(1-s)})$ and $\tau_{P}(m^{*}_{S})=0$, we have
\[
\sqrt{\det A(\tau_{P}(\gamma/N))}
=\sqrt{\det A}(1+O(N^{-(1-s)})).
\]
This combined with Proposition \ref{PRIM} shows the first assertion in Theorem \ref{lCLT}.
Next, we consider the exponent $\delta_{c}(S;\gamma/N)$.
Since $A(\tau)=(\partial \mu_{P})(\tau)$ is bounded from below and since $\tau_{P}=\mu_{P}^{-1}$,
we have
\[
\tau_{P}(x)=\tau_{P}(x)-\tau_{P}(m^{*}_{S})
=A^{-1}(x-m^{*}_{S}) +O(|x-m^{*}_{S}|^{2}).
\]
near $x =m^{*}_{S}$. A Taylor expansion for the function
$f_{\gamma/N}(\tau):=\log k(\tau) -\ispa{\gamma/N,\tau}$ at $\tau =0$ gives
\[
f_{\gamma/N}(\tau)=\log(V(S)) -N^{-1}\ispa{d_{N}(\gamma),\tau} +\ispa{A\tau,\tau}/2 +O(|\tau|^{3}).
\]
Now, noting $A(\tau)=(\partial \mu_{P})(\tau)$, $\tau_{P}=\mu_{P}^{-1}$ and $\tau_{P}(m^{*}_{S})=0$, we have
\[
\tau_{P}(x)=A^{-1}(x-m^{*}_{S}) +O(|x-m^{*}_{S}|^{2}).
\]
These two inequalities with the fact that $\delta_{c}(S;\gamma/N) =f_{\gamma/N}(\tau_{P}(\gamma/N))$ show that
\[
N\delta_{c}(S;\gamma/N)=
N\log(V(S)) -\ispa{A^{-1}d_{N}(\gamma),d_{N}(\gamma)}/(2N) +O(N^{-2}|d_{N}(\gamma)|^{3}).
\]
From this, it is clear that, if $d_{N}(\gamma)=o(N^{s})$ with $0 \leq s \leq 2/3$, then
$O(N^{-2}|d_{N}(\gamma)|^{d})=o(N^{3s-2})$ with $3s-2 \leq 0$, which completes the proof.
\hfill\qedsymbol

\vspace{5pt}

\begin{example}
Let us examine Theorems \ref{lpasympt} and \ref{lCLT} for the one dimensional case.
Let $0 < p$ be a positive integer, and set $S=\{0,1,\ldots,p\}$, so that $P=[0,p]$.
We take $L=\mb{Z} \subset X=\mb{R}$.
Let $c$ be a function on $S$ defined by $c(\beta)={p \choose \beta}$.
In this case, the finite group $\Pi(S)$ is trivial.
The $S$-character $k$, the moment map $\mu_{P}$ and its inverse $\tau_{P}$ are given by
\[
k(\tau)=(1+e^{\tau})^{p},\quad
\mu_{P}(\tau)=\frac{pe^{\tau}}{1+e^{\tau}}, \quad
\tau_{P}(x)=\log (x/(p-x)),\quad x \in (0,p), \tau \in \mb{R}.
\]
Then, the corresponding $S$-lattice paths counting function $\mathcal{P}_{N}^{c}$ is given by
\[
\mathcal{P}_{N}^{c}(\gamma)={Np \choose \gamma}=\frac{(Np)!}{\gamma!(Np-\gamma)!},
\quad \gamma \in [0,Np].
\]
Let $0 < \alpha < p$ be an integer. Then, we have
$\delta (S,\alpha)=p\log p -\alpha \log \alpha -(p-\alpha) \log (p-\alpha)$
and $\det A(S,\alpha) = \frac{\alpha(p-\alpha)}{p}$.
Thus, Theorem \ref{lpasympt} tells us that
\[
\mathcal{P}_{N}^{c}(N\alpha)=\frac{1}{\sqrt{2\pi N}}
\frac{p^{Np+1/2}}{\alpha^{N\alpha +1/2}(p-\alpha)^{N(p-\alpha)+1/2}}(1+O(N^{-1})),
\]
which we can also deduce from Stirling's formula.
Next, a direct computation shows
\[
\delta_{c}(S,x)=\log
\left( \frac{p^{p}}{x^{x}(p-x)^{p-x}}
\right),\quad
A_{c}(S,x)=\frac{x(p-x)}{p},
\quad m^{*}_{S}=\frac{p}{2}.
\]
Therefore, if $\gamma =Np/2 +d_{N}(\gamma)$ is a positive integer
with $d_{N}(\gamma)=o(N^s)$, $0 \leq s \leq 2/3$,
then, by Theorem \ref{lCLT}, we have
\[
\mathcal{P}_{N}^{c}(\gamma)={Np \choose \gamma}
\sim \frac{2^{N+1}}{\sqrt{2\pi pN}}
e^{-\frac{2}{pN}(\gamma -pN/2)^{2}},
\]
which is the fact discussed in Introduction.
\end{example}

\section{Application to multiplicities of group representations}
\label{GREP}

In this section, we shall prove Theorems \ref{SIMPLE},
\ref{repCLT}, \ref{BIANE}, \ref{mwasympt}, \ref{mwasympt2},
\ref{IMULT} as applications of Theorems \ref{lCLT} and
\ref{lpasympt}. As in the introduction, let $G$ be a compact
connected Lie group, and we fix a maximal torus $T$ in $G$. For
any irreducible representation $V_{\lambda}$ of $G$ with highest
weight $\lambda$, the multiplicity of a weight $\nu$ in the
$N$-th tensor power $V_{\lambda}^{\otimes N}$ is denoted by
$m_{N}(\lambda;\nu)$. Similarly, the multiplicity of an
irreducible summand $V_{\nu}$ in $V_{\lambda}^{\otimes N}$ with
the highest weight $\nu$ is denoted by $a_{N}(\lambda;\nu)$.

\subsection{Relation between number of lattice paths and multiplicities}

First of all, we shall explain the relations between the weighted
number of lattice paths discussed in Section \ref{partition} and
the multiplicities $m_{N}$ and $a_{N}$ in group representations.
The main results are Propositions \ref{MequalsP} and
\ref{Mofirred}. In this subsection, we prepare lemmas and
propositions.

Let $\mf{g}$ and $\mf{t}$ be the Lie algebras of $G$ and $T$ respectively.
We fix an inner product $\ispa{\cdot,\cdot}$ on $\mf{g}$ invariant under the adjoint action,
which determines an inner product on $\mf{t}$ invariant under the Weyl group $W$.
In case where $G$ is semisimple, we use the negative Killing form as a fixed inner product.
We sometimes identify
the spaces $\mf{g}$ and $\mf{t}$ with their duals $\mf{g}^{*}$ and $\mf{t}^{*}$, respectively,
by the fixed inner product.
Let $I \subset \mf{t}$ be the integral lattice, {\it i.e.,} $I=\exp^{-1}(1)$, and let
$I^{*} \subset \mf{t}^{*}$ be its dual lattice, {\it i.e.,} the lattice of weights.
We fix an (open) dual Weyl chamber $C$ in $\mf{t}^{*}$.
Let $\Phi$ and $\Phi_{+}$ denote, respectively, the sets of the roots and the positive roots,
respectively.
Let $B \subset \Phi_{+}$ be the set of the simple roots, so that $f \in C$ if and only if
$\ispa{f,\alpha} >0$ for all $\alpha \in B$.
Let $X^{*}$ be the linear span of the simple roots in $\mf{t}^{*}$, and let $X=X^{**}$ be its dual space.
The vector space $X$ is regarded as a subspace in $\mf{t}$ by using the fixed inner product.
Since the simple roots are linearly independent, they form a basis of the vector space $X^{*}$.
Thus we have $\dim X^{*} =\# B =:m$.
The subspace $X \subset \mf{t}$ is spanned by the inverse roots
$\alpha^{*}:=2\kappa^{-1}(\alpha)/\ispa{\alpha,\alpha}$, where $\kappa:\mf{t} \to \mf{t}^{*}$
is an isomorphism induced by the fixed $W$-invariant inner product $\ispa{\cdot,\cdot}$.
We also note that all the roots is in $X^{*}$.

Each dominant weight $\lambda \in \ol{C} \cap I^{*}$ corresponds
to an irreducible unitary representation $V_{\lambda}$.
We define the finite set $M_{\lambda} \subset I^{*}$ by the support of the multiplicity function:
\[
M_{\lambda}:=
\{\mu \in I^{*}\,;\,m_{1}(\lambda;\mu) \neq 0\},
\]
where $m_{1}(\lambda;\mu)$ denotes the multiplicity of the weight $\mu$ in $V_{\lambda}$.
Note that the convex hull $Q(\lambda)$ of the $W$-orbit of $\lambda$
coincides with the convex hull of $M_{\lambda}$.
The dimension of the polytope $Q(\lambda)$ might be less than that of $\mf{t}$.
However, as we shall see soon, the polytope $Q(\lambda)$ is contained in the affine subspace
$X^{*} +\lambda$ in $\mf{t}^{*}$.
Thus, the interior $Q(\lambda)^{o}$ of $Q(\lambda)$ means, in the following,
the interior of $Q(\lambda)$ considered as a polytope in the above affine subspace.
If $G$ is semisimple, then clearly $X^{*}=\mf{t}^{*}$, and hence we can use the polytope $Q(\lambda)$
as the polytope $P$ in Section \ref{partition}.
However, in general, the finite set $M_{\lambda} \subset I^{*}$ of all the weights in $V_{\lambda}$
is not in the subspace $X^{*}$. Thus, we have to modify it. Namely, we set
\[
S_{\lambda}=\{\mu -\lambda\,;\,\mu \in M_{\lambda}\}.
\]

\begin{lem}
We set $D(S_{\lambda})=\{\beta-\beta^{\prime}\,;\,\beta,\beta^{\prime} \in S_{\lambda}\}$.
If $\lambda \in C \cap I^{*}$,
then we have
\[
\lspan_{\mb{R}}D(S_{\lambda})=X^{*},
\]
where the subspace $X^{*} \subset \mf{t}^{*}$ is, as above, the linear span of the simple roots.
\label{subspace}
\end{lem}

\begin{rem}
It should be noted that we denote by $C$ the {\it open} Weyl chamber.
If $\lambda \in \ol{C}$ is contained in a wall,
the linear span $\lspan_{\mb{R}}D(S_{\lambda})$ will be a proper subspace of $X^{*}$.
In fact, in the case where $G=U(2)$, the Weyl group is the symmetric group of order $2!=2$,
and the Weyl chamber is a half-plane in a two dimensional vector space.
Thus, if $\lambda$ is in the wall, which is the unique wall defined by the orthogonal complement of the
(unique) positive root, then it is stable under the Weyl group action.
Thus, the corresponding set $M_{\lambda}$ consists of the single point $\lambda$,
and the linear span $\lspan_{\mb{R}}D(S_{\lambda})$ is the trivial subspace $\{0\}$.
\end{rem}

\begin{proof}
The Weyl group is generated by the reflections $s_{\alpha}$, $\alpha \in B$
with respect to the walls
$H_{\alpha}=\ker \alpha \subset \mf{t} \cong \mf{t}^{*}$, $\alpha \in B$.
Here, we have identified $\mf{t}$ with $\mf{t}^{*}$ by the fixed $W$-invariant inner product.
The reflection $s_{\alpha}$ with respect to the wall $H_{\alpha}$ is given by
\[
s_{\alpha}(f)=f-f(\alpha^{*})\alpha,\quad f \in \mf{t}^{*},
\]
where $\alpha^{*} \in \mf{t}$ is the inverse root of $\alpha$.
Therefore, we have $s_{\alpha}\mu-\mu \in X^{*}$ for every $\mu \in \mf{t}^{*}$ and $\alpha \in B$.
For general $w \in W$, we can write it as a product $w=s_{\alpha_{1}} \cdots s_{\alpha_{n}}$ of reflections.
Then
\[
w\mu -\mu =\sum_{j=1}^{n}(s_{\alpha_{j}}w_{j+1}\mu-w_{j+1}\mu),\quad \mu \in \mf{t}^{*},
\]
where we set $w_{j}=s_{\alpha_{j}}\cdots s_{\alpha_{n}}$ and $w_{n+1}=1$.
The multiplicities of weights in an irreducible representation is invariant under the Weyl group.
Thus, $w\mu \in M_{\lambda}$
for all $\mu \in M_{\lambda}$ and $w \in W$. Therefore, the above expression shows that
$w\mu -\mu \in X^{*}$ for all $\mu \in M_{\lambda}$.
The dominant weight $\lambda$ is in $M_{\lambda}$ with multiplicity $m_{1}(\lambda;w\lambda)=1$, $w \in W$.
Other weights in the representation is lower than $\lambda$.
Thus, we may write
\[
\mu=\sum_{w \in W}c_{w}(\mu)w\lambda,\quad 0 \leq c_{w}(\mu),\quad \sum_{w \in W}c_{w}(\mu)=1,\quad
\mu \in M_{\lambda}.
\]
From this expression, we have $\mu-\lambda \in X^{*}$ for all $\mu \in M_{\lambda}$,
and hence $\lspan_{\mb{R}}D(S_{\lambda}) \subset X^{*}$.
As above, we have $\lambda(\alpha^{*})\alpha =\lambda -s_{\alpha}\lambda$.
But, we have assumed that $\lambda \in C$, the interior of the closed positive Weyl chamber.
Therefore $\ispa{\lambda,\alpha} >0$ for all simple root $\alpha$, and hence $\lambda (\alpha^{*}) >0$,
which implies $\lspan_{\mb{R}}D(S_{\lambda}) = X^{*}$,
completing the proof.
\end{proof}

We consider the lattice $L^{*}=X^{*} \cap I^{*}$ of weights in $X^{*}$
as a fixed lattice in $X^{*}$, as in Section \ref{partition}.
In Section \ref{partition}, the lattice $L(S)^{*}$ spanned by $D(S)$ played a role.
In our case, the lattice $L(S_{\lambda})^{*}$ spanned by $D(S_{\lambda})$ does not depend on $\lambda$
for generic $\lambda$ as follows.

\begin{lem}
Let $\Lambda^{*} \subset X^{*}$ be the lattice spanned by the roots over $\mb{Z}$.
Assume that the dominant weight $\lambda$ is in the open Weyl chamber $C$.
Then we have
\[
L(S_{\lambda})^{*}:=\lspan_{\mb{Z}}(D(S_{\lambda})) =\Lambda^{*}.
\]
\label{ROOTlatt}
\end{lem}

\begin{proof}
First, we note that the irreducible representation $V_{\lambda}$ of $G$ is also irreducible as
a representation of the Lie algebra $\mf{g}$ of $G$ and hence its complexification $\mf{g}^{\mb{C}}$.
Let $\mf{g}_{\alpha}$ denote the root space with the root $\alpha$ so that
\[
\mf{g}^{\mb{C}}=\mf{t}^{\mb{C}} \oplus \bigoplus_{\alpha \in \Phi}\mf{g}_{\alpha}.
\]
We denote by $V_{\lambda}(\mu)$ the weight space in $V_{\lambda}$ for the weight $\mu$.
If $\mu \not \in M_{\lambda}$ then $V_{\lambda}(\mu) =0$.
Consider the subspace $W$ of $V_{\lambda}$ defined by
\[
W=\bigoplus_{\gamma \in \Lambda^{*}}V_{\lambda}(\lambda +\gamma) \subset V_{\lambda}.
\]
It is well-known (\cite{FH}) that the root space $\mf{g}_{\alpha}$ maps $V_{\lambda}(\mu)$ to
$V_{\lambda}(\mu +\alpha)$, and $\mf{t}^{\mb{C}}$ maps $V_{\lambda}(\mu)$ onto itself.
Thus, by the decomposition of $\mf{g}^{\mb{C}}$ above,
the subspace $W$, which contains the one-dimensional subspace $V_{\lambda}(\lambda)$,
is invariant under $\mf{g}^{\mb{C}}$.
Hence we have $W=V_{\lambda}$ by the irreducibility.
Therefore, we have $\mu -\mu' \in \Lambda^{*}$ for every $\mu,\mu' \in M_{\lambda}$, which implies
$L(S_{\lambda})^{*} \subset \Lambda^{*}$.
This holds for arbitrary dominant weight $\lambda \in \overline{C}$.
Now, we assume that $\lambda \in C$.
This implies that the integer $\lambda (\alpha^{*})$ is strictly positive for every simple root $\alpha$.
It is also well-known that the string of weights of the form
\[
\lambda,\ \lambda -\alpha,\ \ldots,\ s_{\alpha}\lambda =\lambda -\lambda(\alpha^{*})\alpha
\]
is contained in $M_{\lambda}$. In particular, we have $\lambda -\alpha \in M_{\lambda}$.
This shows that $\alpha \in L(S_{\lambda})^{*}$ for every simple root $\alpha$.
Since every root can be expressed as a linear combination of the simple roots with
integer coefficients, we have $\Lambda^{*} \subset L(S_{\lambda})^{*}$, and hence the proof.
\end{proof}

By Lemma \ref{subspace}, the finite set $S_{\lambda}$ is a subset in $L^{*}$.
Let $P_{\lambda} \subset X^{*}$ be the convex hull of the finite set $S_{\lambda}$.
The relation of the polytopes $Q(\lambda)$ and $P_{\lambda}$ is
\[
P_{\lambda}=Q(\lambda)-\lambda \subset X^{*}.
\]
The polytope $P_{\lambda}$ contains the origin in $X^{*}$ as a vertex.
Finally, we define the weight function $c_{\lambda}$ on $S_{\lambda}$ by
\[
c_{\lambda}(\beta):=m_{1}(\lambda;\mu),\quad
\beta =\mu-\lambda \in S_{\lambda},
\]
which is, of course, a strictly positive function on $S_{\lambda}$.
Thus, we get the data, $X^{*}$, $L^{*}$, $S_{\lambda}$, $c_{\lambda}$ exactly as in Section \ref{partition}.
Furthermore, we have the following.

\begin{prop}
Let $\mathcal{P}_{N}^{c_{\lambda}}(\gamma)$, $\gamma \in L^{*}$
be the lattice paths counting function in $L^{*}$ with the weight function $c_{\lambda}$
and the set of the allowed steps $S_{\lambda}$.
Then we have
\[
m_{N}(\lambda;\mu)=\mathcal{P}_{N}^{c_{\lambda}}(\mu-N\lambda)
\]
for every $\mu \in NQ(\lambda)$.
\label{MequalsP}
\end{prop}

\begin{proof}
Let $\chi_{\lambda}$ be the character of $V_{\lambda}$, which is considered as a function on $\mf{t}$.
The character $\chi_{\lambda}$ is given explicitly by
\begin{equation}
\chi_{\lambda}(\varphi)=\sum_{\mu \in M_{\lambda}}m_{1}(\lambda;\mu)e^{2\pi i\ispa{\mu,\varphi}},
\quad \varphi \in \mf{t}.
\label{charact1}
\end{equation}
The character of the tensor power $V_{\lambda}^{\otimes N}$ is the $N$-th power $\chi_{\lambda}^{N}$
of the character $\chi_{\lambda}$. Since the multiplicity $m_{N}(\lambda;\mu)$ is the coefficients
of $e^{2\pi i\ispa{\mu,\varphi}}$ in $\chi_{\lambda}^{N}$, we have
\[
m_{N}(\lambda;\mu)=
\sum_{\mu_{1},\ldots,\mu_{N} \in M_{\lambda},\,
\mu=\mu_{1}+\cdots +\mu_{N}}
m_{1}(\lambda,\mu_{1})\cdots m_{1}(\lambda,\mu_{N}).
\]
This shows that $m_{N}(\lambda;\mu)=0$ if $\mu \not \in NQ(\lambda)$.
On the other hand, consider, as in Section \ref{partition}, the
weighted polytope character:
\begin{equation}
k(w)=\sum_{\beta \in S_{\lambda}}
c_{\lambda}(\beta)e^{\ispa{\beta,w}},\quad
w \in X^{\mb{C}}.
\label{pcharact1}
\end{equation}
Then, the lattice paths counting function $\mathcal{P}_{N}^{c_{\lambda}}(\gamma)$ for $\gamma \in L^{*}$
is the coefficient of $e^{\ispa{\gamma,w}}$ in $k(w)^{N}$.
By the definition of the finite set $S_{\lambda}$, we can rewrite the function $k(i\varphi)$
for $\varphi \in X$ as
\begin{equation}
k(i\varphi)=e^{-i\ispa{\lambda,\varphi}}\chi_{\lambda}(\varphi/2\pi),\quad
\varphi \in X (\subset \mf{t}).
\label{CHvsWCH}
\end{equation}
Thus, the coefficient $\mathcal{P}_{N}^{c_{\lambda}}(\mu -N\lambda)$ of $e^{i\ispa{\mu-N\lambda,\varphi}}$
in $k(i\varphi)^{N}$ coincides with $m_{N}(\lambda;\mu)$, concluding the assertion.
\end{proof}

Next, we discuss the multiplicities of irreducible subrepresentations
in the tensor power $V_{\lambda}^{\otimes N}$.
Our strategy to prove Theorem \ref{IMULT} is based on the following alternating sum formula.

\begin{prop}
We fix a dominant weight $\lambda \in C \cap I^{*}$.
Let $\rho$ be half the sum of the positive roots: $\rho =\frac{1}{2}\sum_{\alpha \in \Phi_{+}}\alpha$.
Then we have
\[
\begin{split}
a_{N}(\lambda;\mu)&=\sum_{w \in W}\sgn (w)m_{N}(\lambda;\mu +\rho -w\rho) \\
&=\sum_{w \in W}\sgn (w)\mathcal{P}_{N}^{c_{\lambda}}(\mu -N\lambda +\rho -w\rho),
\quad \mu \in \ol{C}\cap I^{*},
\end{split}
\]
where the weighted lattice paths counting function $\mathcal{P}_{N}^{c_{\lambda}}$ with
the weight function $c_{\lambda}$ and the set of the allowed steps $S_{\lambda}$ in $L^{*}$.
\label{Mofirred}
\end{prop}

\begin{proof}
The second equality follows from Proposition \ref{MequalsP}.
Although the first equality is a special case of the expression $(8)$ in \cite{GM},
we give a proof for completeness.
Consider the character $\chi_{\lambda}^{N}$ of $V_{\lambda}^{\otimes N}$, which has the following
expression:
\begin{equation}
\chi_{\lambda}^{N}=\sum_{\mu \in \ol{C} \cap I^{*}}
a_{N}(\lambda;\mu)\chi_{\mu},
\label{ch111}
\end{equation}
where, $\chi_{\mu}$ is the character of an irreducible representation with the highest weight $\mu$.
By the Weyl character formula, we have
\[
\Delta \chi_{\mu}=\sum_{w \in W}\sgn(w)e^{2\pi iw(\mu +\rho)},
\]
where $\Delta$ is the Weyl denominator $\Delta=\sum_{w \in W}\sgn (w)e^{2\pi i w\rho}$.
Multiplying \eqref{ch111} by the Weyl denominator $\Delta$, we have
\begin{equation}
\Delta \chi_{\lambda}^{N}=
\sum_{\mu \in \ol{C}\cap I^{*},\,w \in W}
\sgn(w)a_{N}(\lambda;\mu)e^{2\pi iw(\mu +\rho)},
\label{tensch1}
\end{equation}
which tells us that the multiplicity $a_{N}(\lambda;\mu)$ for $\mu \in \ol{C} \cap I^{*}$
is the coefficient of $e^{2\pi i(\mu +\rho)}$ in $\Delta \chi_{\lambda}^{N}$.
But, the character $\chi_{\lambda}^{N}$ has the decomposition into the weights for $T$.
Therefore we also have
\begin{equation}
\Delta \chi_{\lambda}^{N}=
\sum_{\gamma \in I^{*},\,w \in W}\sgn(w)m_{N}(\lambda;\gamma)e^{2\pi i(\gamma +w\rho)}.
\label{tensch2}
\end{equation}
In \eqref{tensch2}, the term $e^{2\pi i(\mu +\rho)}$ appears for $\gamma \in I^{*}$ with
$\gamma =\mu +\rho -w \rho$ for every $w \in W$. (Note that $\rho -w\rho$ is a weight for every $w \in W$.)
Therefore, the coefficient of $e^{2\pi i(\mu +\rho)}$
in \eqref{tensch2} is given by
\[
\sum_{w \in W}\sgn(w)m_{N}(\lambda;\mu +\rho -w\rho),
\]
which proves the assertion.
\end{proof}

Next, we assume that $G$ is semisimple. In this case, we simply use the set $M_{\lambda}$
for the finite set $S$ as in Section \ref{partition}.
Furthermore, we have the following

\begin{lem}
Assume that $G$ is semisimple.
Then, for any dominant weight $\lambda$ in the open Weyl chamber $C$,
the center of mass $Q^{*}(\lambda) \in Q(\lambda)$ of the polytope $Q(\lambda)$ defined by
\eqref{COMrep} is the origin.
\label{COMeq0}
\end{lem}

\begin{proof}
First of all, we prove an integral representation for the center of mass.
Since $G$ is semisimple, we may identify the Lie algebra $\mf{g}$ with its dual space by the
Killing form.
Let $\mathcal{O}=\mathcal{O_{\lambda +\rho}}$ denote the (co)adjoint orbit through
the linear form $\lambda +\rho$, where $\rho$ is half the sum of the positive roots.
Then we claim:
\begin{equation}
\int_{\mathcal{O}}\xi\,dm(\xi)=Q^{*}(\lambda),
\label{COMint}
\end{equation}
where the measure $m$ is the symplectic volume measure normalized so that $m(\mathcal{O})=1$.
By a formula for the Fourier transform of the
coadjoint orbit, we have
\begin{equation}
\int_{\mathcal{O}}e^{i\ispa{\xi,H}}\,dm(\xi)
=(-2\pi i)^{d}\vol (\mathcal{O})^{-1}
\frac{\sum_{w \in W}\sgn (w)e^{i\ispa{w(\lambda+\rho),H}}}
{\prod_{\alpha \in \Phi_{+}}\ispa{\alpha,H}},
\label{DHformula}
\end{equation}
where $d=\# \Phi_{+}$ is half the dimension of the coadjoint orbit $\mathcal{O}=\mathcal{O}_{\lambda +\rho}$,
and $\vol(\mathcal{O})$ is the symplectic volume of the orbit.
See \cite{GS} for the proof.
Thus, by the Weyl character formula, we have
\begin{equation}
\int_{\mathcal{O}}e^{i\ispa{\xi,H}}\,dm(\xi)
=(-2\pi i)^{d}\vol(\mathcal{O})^{-1}
\frac{\Delta(H/2\pi)}{\prod_{\alpha \in \Phi_{+}}\ispa{\alpha,H}}
\chi_{\lambda}(H/2\pi),
\label{DHformula2}
\end{equation}
where $\Delta$ is the Weyl denominator: $\Delta(H)=\sum_{w}e^{2\pi i\ispa{w\rho,H}}$, and $\chi_{\lambda}$
is the character of the irreducible representation $V_{\lambda}$.
It is not hard to show that the volume $\vol(\mathcal{O})$ (in our normalization) is given by
\[
\vol(\mathcal{O})=(2\pi)^{d}\frac{\prod_{\alpha \in \Phi_{+}}\ispa{\lambda +\rho,\rho}}
{\prod_{\alpha \in \Phi_{+}}\ispa{\alpha,\rho}}=(2\pi)^{d}(\dim V_{\lambda}).
\]
This can be deduced, for example by setting $H=t\rho$, $t>0$ in \eqref{DHformula2} and then letting $t \to 0$.
(In our notation, the weights of the isotropy representation of $G$ at $\lambda +\rho$
are $2\pi$ times the positive roots $\alpha$.)
Therefore, we have
\[
\int_{\mathcal{O}}e^{i\ispa{\xi,H}}\,dm(\xi)
=\frac{1}{\dim V_{\lambda}}f(H)\chi_{\lambda}(H/2\pi),
\]
where $\chi_{\lambda}$ is the character of $V_{\lambda}$ and the function $f(H)$ is given by
\begin{equation}
f(H)=(-i)^{d}\frac{\Delta(H/2\pi)}{\prod_{\alpha \in \Phi_{+}}\ispa{\alpha,H}}.
\label{rest111}
\end{equation}
It is well-known that
\[
\Delta(H/2\pi)=i^{d}\prod_{\alpha \in \Phi_{+}}
\ispa{\alpha,H}
\left(
\frac{\sin (\ispa{\alpha,H}/2)}{\ispa{\alpha,H}/2}
\right)
=i^{d}\prod_{\alpha \in \Phi_{+}}
\ispa{\alpha,H}(1+O(\ispa{\alpha,H}^{2})).
\]
Hence, we have $f(0)=1$,
and the derivative of the function $f$ at the origin is zero.
Therefore, differentiating the above formula at the origin, we get
\[
i\int_{\mathcal{O}}\xi\,dm(\xi)
=\frac{(\partial \chi_{\lambda})(0)}{\dim V_{\lambda}},
\]
which shows \eqref{COMint}.
By \eqref{COMint}, it is clear that the center of mass $Q^{*}(\lambda)$ is invariant under the
adjoint action. However, the Lie algebra $\mf{g}$ is assumed to be semisimple, and hence $Q^{*}(\lambda)$
must be the origin.
\end{proof}

\begin{rem}
If the Lie algebra $\mf{g}$ is simple, then the above proof is much more simplified.
In fact, the center of mass $Q^{*}(\lambda)$ is clearly invariant under the Weyl group.
But, if $\mf{g}$ is simple, then the Weyl group acts on $\mf{t}^{*}$ irreducibly (see \cite{FH}).
Thus $Q^{*}(\lambda)$ must be the origin.
\end{rem}

\subsection{Proof of Theorem \ref{SIMPLE}}

First of all, we shall prove Theorem \ref{SIMPLE}.
By using Proposition \ref{MequalsP}, we have
\begin{equation}
m_{\lambda,N}=
\frac{1}{V(S_{\lambda})^{N}}
\sum_{\nu-N\lambda \in NP_{\lambda}}
\mathcal{P}_{N}^{c_{\lambda}}(\nu -N\lambda) \delta_{\nu/N}
=
\frac{1}{V(S_{\lambda})^{N}}\sum_{\gamma \in NP_{\lambda}}
\mathcal{P}_{N}^{c_{\lambda}}(\gamma) \delta_{\gamma/N +\lambda},
\label{muT}
\end{equation}
where the weighted volume of the finite set $S_{\lambda}$ is given by
\[
V(S_{\lambda})
=\dim V_{\lambda}=\sum_{\nu-\lambda \in S_{\lambda}}c_{\lambda}(\nu-\lambda),
\quad c_{\lambda}(\nu -\lambda)=m_{1}(\lambda;\nu).
\]

The probability measure $m_{S_{\lambda},N}$ on $X^{*}$,
discussed in Section \ref{partition}, is given by
\begin{equation}
m_{S_{\lambda},N}=
\frac{1}{V(S_{\lambda})^{N}}\sum_{\gamma \in NP_{\lambda}}
\mathcal{P}_{N}^{c_{\lambda}}(\gamma) \delta_{\gamma/N},
\label{muX}
\end{equation}
which is different from $m_{\lambda,N}$ in the term $\delta_{\gamma/N+\lambda}$ and $\delta_{\gamma/N}$.
Thus, for any compact supported continuous function $f$ on $\mf{t}^{*}$, let $f_{\lambda}$ be the function
obtained by translating $f$ by $\lambda$: $f_{\lambda}(x)=f(x+\lambda)$.
Then, we have
\begin{equation}
\int_{X^{*}}f_{\lambda}(x)\,dm_{S_{\lambda},N}
=\int_{\mf{t}^{*}}f(x)\,dm_{\lambda,N}.
\label{TRANS}
\end{equation}
The point $m_{S_{\lambda}}$ is equal to $Q^{*}(\lambda)-\lambda$, where, as in Introduction,
the point $Q^{*}(\lambda)$ is given in \eqref{COMrep},
and hence, by Proposition \ref{lattMl}, we have $m_{\lambda,N} \to \delta_{Q^{*}(\lambda)}$ weakly as $N \to \infty$.
\hfill\qedsymbol

\vspace{10pt}

\subsection{Proof of Theorem \ref{CRAMER} and Corollary \ref{CRAMERM}}

Next, we shall prove Theorem \ref{CRAMER}.
By Proposition \ref{latLDP} and \eqref{CHvsWCH}, the measures $\{m_{S_{\lambda},N}\}$ satisfies
the large deviation principle with the rate function
\[
I_{S_{\lambda}}(x)=\sup_{\tau \in X}
\{\ispa{x +\lambda,\tau}-\log (\chi_{\lambda}(\tau/2\pi i)/(\dim V_{\lambda}))\}.
\]
As in \eqref{TRANS}, we have
$dm_{\lambda,N}=(\phi_{\lambda})_{*}dm_{S_{\lambda},N}$ with $\phi_{\lambda}(x)=x+\lambda$,
namely $m_{\lambda,N}(B)=m_{S_{\lambda},N}(B-\lambda)$. Thus, the measure
$m_{\lambda,N}$ satisfies the large deviation principle with
the rate function $I_{S_{\lambda}}(x-\lambda)=I_{\lambda}(x)$, where the function
$I_{\lambda}(x)$ is given in \eqref{RATEG}, which proves Theorem \ref{CRAMER}.
\hfill\qedsymbol

\vspace{10pt}

To prove Corollary \ref{CRAMERM}, we need the following lemmas.

\begin{lem}
Let $C_{N}(\lambda) \subset \overline{C}$ be a set of dominant weights defined by
\[
C_{N}(\lambda)=\{\mu \in \overline{C} \cap I^{*}\,;\,a_{N}(\lambda;\mu) \neq 0\}.
\]
Then, for a weight $\nu \in I^{*}$, the alternating sum
\begin{equation}
\sum_{\sigma \in W}\sgn(\sigma) m_{N}(\lambda;\nu +\rho -\sigma \rho)=0
\label{VANISH}
\end{equation}
if and only if $\nu +\rho \not \in W(\mu +\rho)$ for every $\mu \in C_{N}(\lambda)$.
\label{vanishL}
\end{lem}

\begin{proof}
First, note that in \eqref{tensch1}, the terms $w(\mu +\rho)$ with $w \in W$ and $\mu \in \overline{C}$ are
all distinct since $\mu +\rho \in \C$ for every $\mu \in \overline{C}$. Thus, in \eqref{tensch1}, the coefficient
of $e^{2\pi i(\nu +\rho)}$ vanish if and only if $\nu+ \rho \not \in W(\mu +\rho)$ for every $\mu \in C_{N}(\lambda)$.
Then, comparing \eqref{tensch1} with \eqref{tensch2}, the coefficient of $e^{2\pi i(\nu +\rho)}$ in \eqref{tensch2}
is give by the alternating sum in \eqref{VANISH}, proving the lemma.
\end{proof}

\begin{lem}
Let $\rho$ be half the sum of the positive roots. For each $w \in W$,
we define a map $\psi_{w,N}:\mf{t}^{*} \to \mf{t}^{*}$ by $\psi_{w,N}(x)=x-(\rho -w\rho)/N$.
Then we have
\[
\sum_{w \in W}\sgn(w)(\psi_{w,N})_{*}dm_{\lambda,N}|_{\overline{C}}
=\frac{B_{N}(\lambda)}{(\dim V_{\lambda})^{N}}
dM_{\lambda,N},
\]
where $|_{\overline{C}}$ denotes the restriction to the closed Weyl chamber $\overline{C}$.
\label{ALTM}
\end{lem}

\begin{proof}
A direct computation with Lemma \ref{vanishL} shows that
\[
\begin{split}
\sum_{w \in W}\sgn(w)(\psi_{w,N})_{*}m_{\lambda,N}
&=\frac{1}{(\dim V_{\lambda})^{N}}\!
\sum_{\nu \in I^{*};\nu +\rho \in W(C_{N}(\lambda)+\rho)}
\sum_{w \in W}\sgn(w)m_{N}(\lambda;\nu+\rho -w\rho)\delta_{\nu/N}. \\
&=\frac{1}{(\dim V_{\lambda})^{N}}
\sum_{\mu \in C_{N}(\lambda)}\sum_{\sigma,w \in W}
\sgn(w)m_{N}(\lambda;\mu +\rho -\sigma^{-1}w\rho)\delta_{\frac{\sigma(\mu+\rho)-\rho}{N}},
\end{split}
\]
where, for the second line, the invariance of the multiplicity $m_{N}(\lambda;\cdot)$ under the
Weyl group has been used.
Now, we restrict the above functional on the closed Weyl chamber $\overline{C}$.
The point $\frac{\sigma(\mu +\rho)-\rho}{N}$ is in $\overline{C}$ if and only if
$\sigma(\mu +\rho) \in \overline{C}+\rho$ since $\overline{C}$ is a cone.
But, in the sum above, $\mu$ is a dominant weight. Thus, only $\sigma=1$ term is in $\overline{C}$.
Thus, the assertion follows from Proposition \ref{Mofirred}.
\end{proof}

\noindent{{\it Completion of proof of Corollary \ref{CRAMERM}}}\hspace{3pt}
First of all, we shall prove upper bound in the large deviation principle.
Note that any $\mu \in C_{N}(\lambda)$ is of order $O(N)$ uniformly, since it is in
the convex polytope $NQ(\lambda)$. By the Weyl dimension formula, we have
\[
\dim V_{\mu} =O(N^{d}),\quad \mu \in C_{N}(\lambda),
\]
with $d$ the number of the positive roots. Then, again the Weyl dimension formula shows
\[
(\dim V_{\lambda})^{N}=\sum_{\mu \in C_{N}(\lambda)}a_{N}(\lambda;\mu)(\dim V_{\mu})
=B_{N}(\lambda)O(N^{d}).
\]
Let $F \subset \overline{C}$ be a closed set. Then, by Lemma \ref{ALTM},
\[
\frac{1}{N}\log (M_{\lambda,N}(F))
=\frac{1}{N}\log
\left(
\sum_{w}\sgn(w)m_{\lambda,N}(F+(\rho -w\rho)/N)
\right) +O(N^{-1}\log N).
\]
For any positive integer $n>0$, we set
\[
F_{n}:=\{x \in \overline{C}\,;\,\inf_{y \in F}|x-y| \leq 1/n\},
\]
which is of course a closed set in $\overline{C}$.
We choose a constant $a>0$ so that $|\rho -w\rho| \leq a$ for every $w \in W$. Then, clearly
$F +(\rho -w\rho)/N \subset F_{t}$ for $a/N \leq 1/n$.
Hence, for every $n$, we have
\[
\frac{1}{N}\log M_{\lambda,N}(F)
\leq \frac{1}{N}\log m_{\lambda,N}(F_{n}) +O(N^{-1}\log N).
\]
Since the measures $m_{\lambda,N}$ satisfies the large deviation principle, we obtain
\[
\limsup_{N \to \infty}\frac{1}{N}\log M_{\lambda,N}(F)
\leq -\inf_{x \in F_{n}}I_{\lambda}(x),
\]
where the rate function $I_{\lambda}(x)$ is given by \eqref{RATEG}.
Now, we claim that
\begin{equation}
\lim_{n \to \infty}a_{n}  =\inf_{x \in F}I_{\lambda}(x), \quad
a_{n}:=\inf_{x \in F_{n}}I_{\lambda}(x),
\label{supinf}
\end{equation}
which will completes the proof, where the existence of the limit in the left hand side
is shown as follows.
The set $F_{n}$ is decreasing: $F_{n} \supset F_{n+1}$, and
the sequence $\{a_{n}\}$ is non-decreasing.
This sequence is bounded from above by $a:=\inf_{x \in F}I_{\lambda}(x)$ because $F=\cap_{n \leq 1}F_{n}$.
Thus, $a_{\infty}:=\lim_{n \to \infty}a_{n}$ exists.
In particular $a \geq a_{\infty}$.
The rate function $I_{\lambda}(x)$ is lower-semicontinuous,
and is {\it good} in the sense that
its sublevel set $I_{\lambda}^{-1}[0,\alpha]$ is compact for every $\alpha>0$ (see \cite{DZ}).
Thus, the function $I_{\lambda}$ attains its minimum on each closed set.
Let $x_{n} \in F_{n}$ be a point such that $I_{\lambda}(x_{n})=a_{n}$. Note that $x_{n}$ is in the compact set $I_{\lambda}^{-1}[0,a]$,
and hence it has a convergent subsequence. We also denote it by $x_{n}$.
Since $F$ is closed, there exists a point $y_{n} \in F$ such that $\inf_{y \in F}|y-x_{n}|=|y_{n}-x_{n}| \leq 1/n$,
and, as a result, $\{y_{n}\}$ contains a convergent sequence. Therefore, the limit $x:=\lim x_{n}$ is in $F$.
By the lower-semicontinuity, we have
\[
a_{\infty}=\lim_{n \to \infty}I_{\lambda}(x_{n}) \geq I_{\lambda}(x) \geq a=\inf_{x \in F}I_{\lambda}(x) \geq a_{\infty},
\]
which establishes \eqref{supinf}.
\hfill\qedsymbol

\vspace{10pt}

\subsection{Proof of Theorems \ref{mwasympt} and \ref{mwasympt2}}

By Theorem \ref{lpasympt} and
Proposition \ref{MequalsP},
we have an asymptotic estimate of the multiplicity $m_{N}(\lambda;N\nu +f)$
if $\nu_{0} \in Q(\lambda)^{o}$ and $f \in \Lambda^{*}$.
To compute the exponent $\delta_{c_{\lambda}}(S_{\lambda},\nu_{0}-\lambda)$ and the
linear transform $A_{c_{\lambda}}(S_{\lambda},\nu_{0}-\lambda)$ from $X$ to $X^{*}$ in Theorem \ref{lpasympt},
we note that the moment map \eqref{moment} for $S=S_{\lambda}$ is given by
\[
\mu_{P_{\lambda}}:X \ni x \to \mu_{\lambda}(x)-\lambda \in P_{\lambda},
\]
where $\mu_{\lambda}$ is defined in \eqref{Gmomexp}.
Thus, we have $\tau_{\lambda}(\nu_{0})=\tau_{P_{\lambda}}(\nu_{0}-\lambda)$.
From this, we have $\delta_{c_{\lambda}}(S_{\lambda},\nu_{0}-\lambda)=\delta_{\lambda}(\nu_{0})$.
The positivity of the linear transform $A_{c_{\lambda}}(S_{\lambda},\nu_{0}-\lambda)$ from $X$ to $X^{*}$
is proved in Section \ref{partition}. A direct computation by using the definition \eqref{matrix0} shows that
\begin{gather*}
A_{c_{\lambda}}(S_{\lambda},\nu_{0}-\lambda)
=\sum_{\mu \in M_{\lambda}}k_{\mu}(\nu_{0})
(\mu -\lambda) \otimes (\mu -\lambda)
-(\nu_{0} -\lambda) \otimes (\nu_{0} -\lambda), \\
k_{\mu}(\nu_{0}):=\frac{m_{1}(\lambda;\mu)e^{\ispa{\mu,\tau_{\lambda}(\nu_{0})}}}
{\sum_{\mu' \in M_{\lambda}}m_{1}(\lambda;\mu')e^{\ispa{\mu',\tau_{\lambda}(\nu_{0})}}}
\end{gather*}
where, for any $f \in X^{*}$, $f \otimes f :X \to X^{*}$ is defined by $(f \otimes f)x=\ispa{x,f}f$, $x \in X$.
By definition (\eqref{Gmomexp}), we have $\sum_{\mu}k_{\mu}(\nu_{0})\mu=\mu_{\lambda}(\tau_{\lambda}(\nu_{0}))=\nu_{0}$.
From this, it is easy to see that $A_{c_{\lambda}}(S_{\lambda},\nu_{0}-\lambda)$ coincides with the
linear transform $A^{0}_{\lambda}(\nu_{0})$ on $X$. This shows that $A_{\lambda}(\nu_{0})$ is positive definite
as a linear transform from $X$ to $X^{*}$, and it is equal to $A_{c_{\lambda}}(S_{\lambda},\nu_{0}-\lambda)$.
The positivity of the exponent $\delta_{\lambda}(\nu_{0})$ follows from the assumption that
the weight $\nu_{0}$ occurs in $V_{\lambda}$. This completes the proof of
Theorem \ref{mwasympt}. Similarly, Theorem \ref{mwasympt2} is proved by using Proposition \ref{PRIM}.
\hfill\qedsymbol

\vspace{10pt}

\subsection{Proof of Theorem \ref{repCLT}}

Before proving Theorem \ref{repCLT}, we shall state more general result, which
corresponds to Theorem \ref{lCLT}.
\begin{theo}\label{MDgen}
Let $0 \leq s \leq 2/3$. Let $\nu \in NQ(\lambda)$ be a weight of the form
\[
\nu =NQ^{*}(\lambda) +d_{N}(\nu),\quad |d_{N}(\nu)|=o(N^{s}).
\]
Assume that $m_{N}(\lambda;\nu) \neq 0$ for every sufficiently large $N$.
Then, we have
\[
m_{N}(\lambda;\nu)=(2\pi N)^{-m/2}|\Pi(G)|(\dim V_{\lambda})^{N}
\frac{e^{-\ispa{A_{\lambda}^{-1}d_{N}(\nu),d_{N}(\nu)}/(2N)}}
{\sqrt{\det A_{\lambda}}}(1+\varepsilon_{N}),
\]
where
\[
\varepsilon_{N}=
\left\{
\begin{array}{ll}
O(N^{-(1-s)}) & \  \mbox{for}\ 0 \leq s \leq 1/2, \\
o(N^{3s-2}) & \ \mbox{for} \ 1/2 < s \leq 2/3,
\end{array}
\right.
\]
and the positive definite linear transformation $A_{\lambda}:X \to X^{*}$ is given by
\[
A_{\lambda}=A_{\lambda}(Q^{*}(\lambda))=
\frac{1}{\dim V_{\lambda}}\sum_{\mu \in M_{\lambda}}m_{1}(\lambda;\mu)\mu \otimes \mu
-Q^{*}(\lambda) \otimes Q^{*}(\lambda).
\]

\end{theo}

\begin{proof}
This follows from Theorem \ref{lCLT} and Proposition \ref{MequalsP}, and the computations for
the exponent and the matrix by the same method as in the proof of Theorem \ref{mwasympt}.
\end{proof}

\noindent{{\it Completion of Proof of Theorem \ref{repCLT}}}\hspace{3pt}
Assume that $G$ is semisimple. Then, by Lemma \ref{COMeq0}, $Q^{*}(\lambda)=0$.
Thus, $d_{N}(\lambda)$ is $\gamma$ itself. Hence, Theorem \ref{repCLT}
is a direct consequence of Theorem \ref{MDgen}.
\hfill\qedsymbol

\subsection{Proof of Theorems \ref{IMULT} and \ref{BIANE}}

For any $w \in W$, the weight $\rho -w\rho$ is in the root lattice $\Lambda^{*}$.
Therefore, we can apply Theorem \ref{mwasympt} for $f=\rho -w\rho$ and $\nu_{0}=\nu$.
Now, Theorem \ref{IMULT} follows from Proposition \ref{Mofirred}
\hfill\qedsymbol

\subsubsection{Proof of Theorem \ref{BIANE}}

As mentioned in Introduction, our approach to the irreducible
multiplicities based on Proposition \ref{Mofirred} does not seem
to be the most efficient for the central limit region. Our
steepest descent method easily gives the principal term, but the
remainder estimate becomes tricky since one needs to use
cancellations occuring in the alternating sum over the Weyl group.
Hence, we use the method of Biane \cite{B} in this region.
Although it is not new, we include it for the sake of
completeness. We also add some details not in \cite{B}.

\vspace{5pt}

We begin with:

\begin{lem}
Assume that $G$ is semisimple, and assume also that the fixed
dominant weight $\lambda$ is in the open Weyl chamber $C$. For
$N>0$, we set
\[
NM_{\lambda}=\{\mu=\nu_{1} +\cdots+\nu_{N}\,;\,\nu_{j} \in
M_{\lambda},j=1,\ldots,N\}.
\]
Let $\mu$ be a dominant weight such that $\mu \not \in
NM_{\lambda}+\Lambda^{*}$. Then $a_{N}(\lambda;\mu)=0$.
\label{NULLCOND}
\end{lem}

\begin{proof}
Since $G$ is assumed to be semisimple, we may use the polytope
$Q(\lambda)$ as $P$ in Section \ref{partition} and $M_{\lambda}$
as the finite set $S$. Thus, the torus ${\bf T}^{m}$ essentially
coincides with the maximal torus $T$. The finite group
$\Pi(G)$ is isomorphic to the kernel of the surjective
homomorphism $\pi_{\lambda}:{\bf T}^{m} \to T(G):=X/(2\pi
\Lambda)$. We also note that $\Lambda^{*} \subset I_{\lambda}^{*}
\subset L^{*}$, where $L^{*}=I^{*}$ is the full weight lattice,
where $I_{\lambda}^{*}$ is the lattice spanned by $M_{\lambda}$
over $\mb{Z}$.

By the Weyl integration formula (or by using Propositions
\ref{MequalsP}, \ref{Mofirred} and the integral formula
\eqref{integral1}), we have
\begin{equation}
a_{N}(\lambda;\mu)=\frac{(\dim V_{\lambda})^{N}}{(2\pi)^{m}}
\int_{{\bf T}^{m}}e^{-i\ispa{\mu +\rho,\varphi}}
K(\varphi)^{N}J(\varphi)\,d\varphi, \label{INTI}
\end{equation}
where we set $K(\varphi)=\chi_{\lambda}(\varphi/2\pi)/\dim
V_{\lambda}$ and $J(\varphi)=\Delta (\varphi/2\pi)$ being
$\chi_{\lambda}$ the character of $V_{\lambda}$ and $\Delta$ the
Weyl denominator $\Delta (H)=\sum_{w \in W}\sgn (w) e^{2\pi
i\ispa{w\rho,H}}$. As in the proof of Theorem \ref{lpasympt}
(Section \ref{partition}), we use the cut-off function $\chi$
around the origin so that a branch of the logarithm $\log K$
exists on $\supp \chi$. We also use the function
$\chi_{g}=\chi(\varphi-\varphi_{g})$, where $\varphi_{g} \in 2\pi
\Lambda$ is a (fixed) representative of $g \in \ker \pi_{\lambda}
\cong \Pi(G)$, {\it i.e.}, $g =\exp \varphi_{g} \in {\bf T}^{m}$
$\pi_{\lambda}(\exp \varphi_{g})=1$. Then, by Lemma
\ref{equalone}, we have
\[
a_{N}(\lambda;\mu)=\frac{(\dim V_{\lambda})^{N}}{(2\pi)^{m}}
\left( \sum_{g \in \ker \pi_{\lambda}} \int e^{-i\ispa{\mu
+\rho,\varphi}}K(\varphi)^{N}J(\varphi)\chi_{g}(\varphi)\,d\varphi
+O(e^{-cN}) \right)
\]
for some constant $c>0$. Now, we make a change of variable
$\varphi \mapsto \varphi +\varphi_{g}$ for each integral in the
above. Then, we will have the term
\begin{equation}
e^{-i\ispa{\mu +\rho,\varphi_{g}}}h(g)^{N}J(\varphi +\varphi_{g})
=\sum_{w \in W}\sgn (w) [e^{-i\ispa{\mu +\rho,\varphi_{g}}}
h(g)^{N}e^{i\ispa{w\rho,\varphi_{g}}} ]e^{i\ispa{w\rho,\varphi}}
\label{WWW1}
\end{equation}
in the integrand, where $h(g)=e^{i\ispa{\nu,\varphi_{g}}}$ for $g
\in \ker \pi_{\lambda} \cong \Pi(G)$ which does not depends on the
choice of $\nu \in M_{\lambda}$. Note $\rho -w\rho \in
\Lambda^{*}$ for every $w \in W$. Thus $\ispa{\rho
-w\rho,\varphi_{g}}$ is $2\pi$ times an integer. Therefore, the
expression \eqref{WWW1} is equal to
$e^{-i\ispa{\mu,\varphi_{g}}}h(g)^{N}J(\varphi)$. Hence we obtain
\[
a_{N}(\lambda;\mu)= \frac{(\dim V_{\lambda})^{N}}{(2\pi)^{m}}
\left( \sum_{g \in \Pi(G)}e^{-i\ispa{\mu,\varphi_{g}}}h(g)^{N}
\right) \int e^{-i\ispa{\mu
+\rho,\varphi}}K(\varphi)^{N}J(\varphi)\chi(\varphi)\,d\varphi
\]
modulo $O(e^{-cN})$. The map $g \mapsto
e^{-i\ispa{\mu,\varphi_{g}}}h(g)^{N}$ is a character of the finite
abelian group $\Pi(G)$. This character is not a trivial character
of $\Pi(G)$ if and only if
\[
\sum_{g \in \Pi(G)}e^{-i\ispa{\mu,\varphi_{g}}}h(g)^{N}= \sum_{g
\in \Pi(G)}e^{i\ispa{N(\nu_{1}+\cdots+\nu_{N})
-\mu,\varphi_{g}}}=0,
\]
where $\nu_{j} \in M_{\lambda}$, $j=1,\ldots,N$ are weights in $V_{\lambda}$.
The left hand side above does not depend on the choice of $\nu_{j} \in M_{\lambda}$,
$j=1,\ldots,N$. Now, assume that $\mu \not \in NM_{\lambda}
+\Lambda^{*}$, which implies $\mu -N(\nu_{1}+\cdots+\nu_{N})
\not\in \Lambda^{*}$ for every $\nu_{j} \in M_{\lambda}$. Since
$\varphi_{g}$ is arbitrarily fixed representative for each $g \in
\Pi(G)=\Lambda/L$, some of $g \in \Pi(G)$ must satisfy
$\ispa{\mu-N\nu, \varphi_{g}} \not \in 2\pi \mb{Z}$, and hence the
character $e^{-i\ispa{\mu,\varphi_{g}}}h(g)^{N}$ is non-trivial.
Thus, we have $a_{N}(\lambda;\mu)=0$ for $\mu \not \in
NM_{\lambda}+\Lambda^{*}$.
\end{proof}

\begin{rem}
As in the proof above, for $\mu \in NM_{\lambda}+\Lambda^{*}$, one
has
\begin{equation}
a_{N}(\lambda;\mu)= \frac{(\dim
V_{\lambda})^{N}|\Pi(G)|}{(2\pi)^{m}} \int e^{-i\ispa{\mu
+\rho,\varphi}}K(\varphi)^{N}J(\varphi)\chi(\varphi)\,d\varphi
+O(e^{-cN}). \label{MULTINT}
\end{equation}
\end{rem}

We now complete the proof of Theorem \ref{BIANE}. We use the
integral representation \eqref{MULTINT}. For $\mu \in
NM_{\lambda}+\Lambda^{*}$, we have
\begin{gather*}
a_{N}(\lambda;\mu)=\frac{|\Pi(G)|(\dim V_{\lambda})^{N}}{(2\pi)^{m}N^{m/2}}I(N),\\
I(N):=\int e^{-i\ispa{\mu +\rho,\varphi}/N^{1/2}}
K(\varphi/N^{1/2})^{N}J(\varphi/N^{1/2})\chi(\varphi/N^{1/2})\,d\varphi
\end{gather*}
modulo $O(e^{-cN})$.
As in \cite{B}, we set $\kappa (\varphi)=\prod_{\alpha>0}\ispa{\alpha,\varphi}$, which is a
polynomial of degree $d=\# \Phi_{+}$, the number of the positive roots.
Then, it is easy to show that $J(\varphi)=(\frac{i}{N^{1/2}})^{d}\kappa(\varphi)(1+|\varphi|^{2d}O(N^{-1}))$.
Since $|K(\varphi)|^{2}$ is real, and since the first derivative of $K$ at the origin is zero (Lemma \ref{COMeq0}),
we can choose $r>0$ such that
\begin{equation}
|K(\varphi)|^{2} \leq 1 -c\ispa{A_{\lambda}\varphi,\varphi} \leq e^{-c\ispa{A_{\lambda}\varphi,\varphi}},\quad
|\varphi| <r.
\label{RRR1}
\end{equation}
Replacing $\chi$ by a cut-off function whose support is small enough,
we have
\[
\int |K(\varphi/N^{1/2})|^{N}|\kappa(\varphi)||\varphi|^{2d}\chi(\varphi/N^{1/2})\,d\varphi
=O(1),
\]
and hence
\[
I(N)=(i/N^{1/2})^{d}I_{1}(N)(1+O(1/N)),\;\;
I_{1}(N)=\int e^{-i\ispa{\mu +\rho,\varphi}/N^{1/2}}K(\varphi/N^{1/2})\kappa(\varphi)\chi (\varphi/N^{1/2}).
\]
For simplicity, we set $A_{N}(\varphi)=e^{-i\ispa{\mu +\rho,\varphi}/N^{1/2}}\kappa(\varphi)$.
A Taylor expansion of $\log K$ at the origin gives
\begin{equation}
K(\varphi/N^{1/2})^{N}=e^{-\ispa{A_{\lambda}\varphi,\varphi}/2-iT(\varphi)/N^{1/2}}e^{NR_{4}(\varphi/N^{1/2})},
\label{TTT1}
\end{equation}
where $R_{4}(\varphi)=O(|\varphi|^{4})$ locally uniformly.
Concerning this expansion, we write
\begin{equation}
I_{1}(N)=\int A(\varphi)e^{-\ispa{A_{\lambda}\varphi,\varphi}/2-iT(\varphi)/N^{1/2}}\,d\varphi
+\sum_{j=1}^{3}I\!\!I_{j}(N),
\label{III1}
\end{equation}
where we set
\begin{gather*}
I\!\!I_{1}(N)  =
\int A(\varphi)(K(\varphi/N^{1/2})-e^{-\ispa{A_{\lambda}\varphi,\varphi}/2-iT(\varphi)/N^{1/2}})
\chi (\varphi/N^{1/4})\,d\varphi, \\
I\!\!I_{2}(N)  =
\int A(\varphi) K(\varphi/N^{1/2})^{N}(1-\chi (\varphi/N^{1/4}))\chi (\varphi/N^{1/2})\,d\varphi,\\
I\!\!I_{3}(N)  =
-\int A(\varphi) e^{-\ispa{A_{\lambda}\varphi,\varphi}/2-iT(\varphi)/N^{1/2}}
(1-\chi (\varphi/N^{1/4}))\,d\varphi.
\end{gather*}
Here we note that $\chi (\varphi/N^{1/4})\chi (\varphi/N^{1/2})=\chi(\varphi/N^{1/4})$ for sufficiently large $N$.
For the integral $I\!\!I_{1}(N)$, the integrand vanish for $|\varphi| >cN^{1/4}$ for some $c$.
Thus, by \eqref{TTT1}, we have $|e^{NR_{4}(\varphi/N^{1/2})}| =O(1)$, and $NR_{4}(\varphi/N^{1/2})=|\varphi|^{4}O(1/N)$.
Therefore we have $|I\!\!I_{1}(N)|=O(1/N)$. For the integral $I\!\!I_{2}(N)$, $\varphi/N^{1/2}$ is bounded.
Thus, by \eqref{RRR1}, we have
\[
|I\!\!I_{2}(N)| \leq \int_{|\varphi| \geq N^{1/4}}e^{-c\ispa{A_{\lambda}\varphi,\varphi}/2}|\kappa(\varphi)|\,d\varphi
=O(N^{(d+m-1)/4}e^{-cN^{1/2}}).
\]
Similarly, it is easy to see that $I\!\!I_{3}(N)=O(N^{(m-2)/4}e^{-cN^{1/2}})$.
Finally, we consider the first integral in \eqref{III1}, which can be written in the form
\[
\int A(\varphi)e^{-\ispa{A_{\lambda}\varphi,\varphi}/2-iT(\varphi)/N^{1/2}}\,d\varphi
=\int e^{-i\ispa{\mu +\rho,\varphi}/N^{1/2}}\kappa (\varphi)e^{-\ispa{A_{\lambda}\varphi,\varphi}/2}
\,d\varphi(1+O(1/N^{1/2})).
\]
By using the homogeneity of the polynomial $\kappa$ of degree $d$, It is easy to see that
\[
\int e^{-i\ispa{\mu +\rho,\varphi}/N^{1/2}}\kappa (\varphi)e^{-\ispa{A_{\lambda}\varphi,\varphi}/2}\,d\varphi
=\frac{i^{d}(2\pi)^{m/2}}{\sqrt{\det A_{\lambda}}}
\kappa (\partial)(e^{-\ispa{A_{\lambda}^{-1}\varphi,\varphi}/2})((\mu +\rho)/N^{1/2}).
\]
As in \cite{B}, by using the fact that the polynomial $\kappa$
is alternating with respect to the $W$-action,
it is not hard to see that
\[
\kappa(\partial)(e^{-\ispa{A_{\lambda}^{-1}\varphi,\varphi}/2})=
(-1)^{d}\kappa (A_{\lambda}^{-1}\varphi).
\]
Therefore, we have
\[
a_{N}(\lambda;\mu)
=\frac{|\Pi(G)|(\dim V_{\lambda})^{N}}
{(2\pi)^{m/2}N^{d+m/2}\sqrt{\det A_{\lambda}}}
\kappa (A_{\lambda}^{-1}(\mu +\rho))(1+O(1/N^{1/2})).
\]
Note that the inner product $\ispa{A_{\lambda}^{-1}x,y}$ is invariant under the action of the Weyl group.
Therefore, by the Weyl dimension formula, we have
\[
\kappa (A_{\lambda}^{-1}(\mu+\rho))=
(\dim V_{\mu})\prod_{\alpha >0}\ispa{A_{\lambda}^{-1}\rho,\alpha},
\]
which concludes the assertion.
\hfill\qedsymbol

\section{Example: $G=U(2)$}
\label{unitary}

In the previous sections, we have obtained the asymptotics of the multiplicities of weights
and irreducibles in high tensor power $V_{\lambda}^{N}$ of a fixed irreducible
representation $V_{\lambda}$.

The leading term of our asymptotic formula are described by
the constant $\delta_{\lambda}(\nu)$ and the determinant
$\det A_{\lambda}(\nu)$ of the matrix $A_{\lambda}(\nu)$.
In general, it seems somewhat difficult to calculate them explicitly.
The most subtle point is the inverse of the ``moment map'' $\tau_{\lambda}(\nu) \in X$.
Furthermore, in Theorem \ref{IMULT}, the term of the Weyl denominator might vanish.
The aim of this section is to discuss them for the group $G=U(2)$.

Roughly speaking, for $G=U(2)$, the corresponding lattice paths model
is {\it Example} in Section \ref{partition} with the weight function $c \equiv 1$.
(But for general $G=U(m+1)$, it is not identically $1$.)

To begin with, we recall some of facts about representation theory for $G=U(m+1)$.
Let $T \subset U(m+1)$ ($m \geq 1$) be the maximal torus of all diagonal matrices in
the unitary group $U(m+1)$.
The Lie algebra $\mf{t}$ of $T$ consists of all diagonal matrices with pure imaginary entries.
We identify $\mf{t}$ with $\mb{R}^{m+1}$ by
$(x_{1},\ldots,x_{m+1}) \mapsto 2\pi i {\rm diag}(x_{1},\ldots,x_{m+1})$.
Let $e_{j}$ ($j=1,\ldots,m+1$) be the standard basis for $\mb{R}^{m+1}$, and let $e_{j}^{*}$ be the dual basis.
The Weyl group $W$ is the symmetric group $\mf{S}_{m+1}$ of order $(m+1)!$.
We use the usual Euclidean inner product to identify $\mf{t} \cong \mb{R}^{m+1}$ with its dual.
The integer lattice and the lattice of weights are identified with $\mb{Z}^{m+1}$.
We choose the positive open Weyl chamber $C$ given by
\[
C=\{\gamma =(\gamma_{1},\ldots,\gamma_{m+1})\,;\,\gamma_{1} > \cdots > \gamma_{m+1}\}.
\]
The roots of $(G,T)$ are $\alpha_{i,j}:=e_{i}^{*}-e_{j}^{*}$, $i \neq j$; the positive roots;
$\alpha_{i,j}$, $i <j$, and the simple roots; $\alpha_{j}:=\alpha_{j,j+1}$, $j=1,\ldots,m$.
The subspace $X^{*} \subset \mf{t}^{*} \cong \mb{R}^{m+1}$ spanned by the simple roots is identified with
\[
X \cong X^{*} = \{(x_{1},\ldots,x_{m+1}) \subset \mb{R}^{m+1}\,;\,\sum x_{j} =0\},
\]
which is identified with the Lie algebra of $T \cap SU(m+1)$.
Half the sum of the positive roots $\rho$ is given by
\begin{equation}
\rho :=\frac{1}{2}\sum_{1 \leq i <j\leq m+1}\alpha_{i,j}
=\frac{1}{2}\sum_{j=1}^{m}(m+2-2j)e_{j}^{*}.
\label{hsumunitary}
\end{equation}
The alternating sum $A(\gamma)$ for the functional $\gamma \in \mf{t}^{*}$ is a function
on $\mf{t}$ given by
\[
A(\gamma)(\varphi):=\sum_{w \in \mf{S}_{m+1}}\sgn(w)e^{2\pi i\ispa{w\gamma,\varphi}},\quad
\varphi \in \mf{t} \cong \mb{R}^{m+1}.
\]
Then the Weyl character formula states that,
for a dominant weight $\lambda \in \ol{C} \cap \mb{Z}^{m+1}$,
the character $\chi_{\lambda}$ for the irreducible representation $V_{\lambda}$ corresponding to $\lambda$
is given by
\[
\chi_{\lambda}(\varphi)=
\frac{A(\lambda +\rho)(\varphi)}{\Delta (\varphi)},\quad
\varphi \in \mf{t},
\]
where $\Delta$ is the Weyl denominator $\Delta =A(\rho)$.
In the case where $G=U(m+1)$,
one can compute the alternating sum $A(\gamma)$ from the definition,
and, as a result,
the character $\chi_{\lambda}$ is given by the {\it Schur polynomial} $s_{\zeta_{\lambda}}$
for the partition $\zeta_{\lambda}=(\lambda_{1}-\lambda_{m+1},\ldots,\lambda_{m}-\lambda_{m+1},0)$:
\begin{gather*}
\chi_{\lambda}(\varphi)=
(\xi_{1}\cdots \xi_{m+1})^{\lambda_{m+1}}
s_{\zeta_{\lambda}}(\xi_{1}(\varphi),\ldots,\xi_{m+1}(\varphi)),\\
s_{\zeta_{\lambda}}:=\frac{\det (\xi_{i}(\varphi)^{(\lambda_{j}-\lambda_{m+1})+m+1-j})}
{\det (\xi_{i}(\varphi)^{m+1-j})},
\quad \xi_{j}:=e^{2\pi i e_{j}^{*}},
\end{gather*}
where the denominator in the above is Vandermond's determinant (difference product):
\[
D(\xi_{1},\ldots,\xi_{m}):=\prod_{1 \leq i<j \leq m+1}(\xi_{i}-\xi_{j}).
\]
If $\lambda_{m+1} \geq 0$, then the above is just the Schur polynomial $s_{\lambda}$
with the partition $\lambda$.

Now we fix a dominant weight $\lambda \in C \cap \mb{Z}^{m+1}$.
For simplicity, we assume that $\lambda_{m+1} \geq 0$ so that the character $\chi_{\lambda}$
is precisely the Schur polynomial $s_{\lambda}$.

It is well-known (see \cite{FH}) that the multiplicity $m_{1}(\lambda;\mu)$ of a partition
$\mu$ (which is equivalent to say that $\mu$ is a dominant weight with non-negative entries)
is given by the {\it Kostka number} $K_{\lambda \mu}$
which is the coefficients in the Schur polynomial $s_{\lambda}$ of the
symmetric sum of the monomials corresponding to $\mu$.
It is also well-known (\cite{FH}) that $K_{\lambda\mu} \neq 0$ if and only if
the partition $\mu$ satisfies
\begin{equation}
\sum_{j=1}^{i}\mu_{j} \leq \sum_{j=1}^{i}\lambda_{j},\quad
i=1,\ldots,m,
\label{kostka}
\end{equation}
and $\sum_{j=1}^{m+1}\mu_{j}=\sum_{j=1}^{m+1}\lambda_{j}$.
(The last condition is necessary, since the weights in $V_{\lambda}$
is in the convex hull of the $W$-orbit of $\lambda$.)

We note that the relation between our weighted character
function $k$ and the character $\chi_{\lambda}$ is expressed as
\begin{equation}
k(\tau)=e^{-\ispa{\lambda,\tau}}\chi_{\lambda}(\tau/2\pi i)=
e^{-\ispa{\lambda,\tau}}s_{\lambda}(e^{\tau_{1}},\ldots,e^{\tau_{m}}),\quad
\tau=(\tau_{1},\ldots,\tau_{m}) \in X (\subset \mf{t}).
\label{schurch}
\end{equation}
Note that, in the above, the character $\chi_{\lambda}$
is extended to the complexification $\mf{t}^{\mb{C}}$.
In particular, we have
\begin{equation}
\log k(\tau)- \ispa{\nu -\lambda,\tau}=
\log s_{\lambda}(e^{\tau}) -\ispa{\nu,\tau},\quad
\tau \in X.
\label{logschur}
\end{equation}
Therefore, as in \eqref{Const}, \eqref{aux11}, the constant $\delta_{\lambda}(\nu)$ is given by
\begin{equation}
\delta_{\lambda}(\nu)=\log s_{\lambda}(e^{\tau_{\lambda}(\nu)})-\ispa{\nu,\tau_{\lambda}(\nu)}.
\label{deltau2}
\end{equation}

Now, consider the case where $m=1$.
We take $\lambda =(\lambda_1,\lambda_2) \in C \cap \mb{Z}^{2}$, $\lambda_1 >\lambda_2 \geq 0$.
We set $n_{\lambda}=\lambda_1 -\lambda_2 >0$. Then, the Schur polynomial $s_{\lambda}(\xi_1,\xi_2)$
in two variables corresponding to the partition $\lambda$ is given by
\begin{equation}
s_{\lambda}(\xi_1,\xi_2)=
\frac{\xi_{1}^{\lambda_{1}+1} \xi_{2}^{\lambda_{2}}-\xi_{1}^{\lambda_{2}}\xi_{2}^{\lambda_{1}+1}}
{\xi_{1}-\xi_{2}}
=\sum_{j=0}^{n_{\lambda}}
\xi_{1}^{\lambda_{1}-j}\xi_{2}^{\lambda_{2}+j}.
\label{schuru2}
\end{equation}
Therefore, the weights in the irreducible representation
$V_{\lambda}$ are of the form:
\begin{equation}
\nu_{j}:=\lambda -j\alpha,\quad j=0,\ldots,n_{\lambda},
\label{weightu2}
\end{equation}
where $\alpha$ is the unique positive (simple) root $\alpha=(1,-1)$.
All these weights have multiplicity one: $m_{1}(\lambda;\nu_{j})=1$.
Therefore, the multiplicity for the high tensor power $V_{\lambda}^{\otimes N}$
is given by (see Proposition \ref{MequalsP})
\[
m_{N}(\lambda;\mu)=\# \{(j_{1},\ldots,j_{N}) \,;\,
0 \leq j_{k} \leq n_{\lambda}, \mu =N\lambda -(j_{1} +\cdots j_{N})\alpha\}.
\]
The polytope $P_{\lambda}$ is given by
\[
P_{\lambda}=\{\tau\alpha \in X^{*}\,;\,-n_{\lambda} \leq \tau \leq 0\}.
\]
Thus, we have the following

\begin{lem}
For every $j=0, \ldots, n_{\lambda}$, $\nu_{j}$ is a weight
in the interior of $Q(\lambda)=P_{\lambda}+\lambda$
if and only if $1 \leq j \leq n_{\lambda}-1$.
Furthermore, $\nu_{j}$ is a dominant weight in the interior of $Q(\lambda)$
if and only if $1 \leq j \leq \frac{n_{\lambda}}{2}$.
\end{lem}

Next, we shall calculate the moment map $\mu_{P_{\lambda}}:X \to P_{\lambda}$ defined in \eqref{moment}.

\begin{lem}
We identify $X^{*}$ with $\mb{R}$ through the identification $\mb{R} \ni \tau \mapsto \tau \alpha \in X^{*}$.
We set $h(\tau) =k(\tau \alpha)$.
Then the moment map $\mu_{P_{\lambda}}$ is given by
\begin{equation}
\mu_{P_{\lambda}}(\tau \alpha) =f(\tau)\alpha,\quad
f(\tau)=\frac{h'(\tau)}{2h(\tau)}.
\label{momentu2}
\end{equation}
The functions $h(\tau)$ and $f(\tau)$ are given explicitly by
\begin{gather*}
h(\tau)=e^{-n_{\lambda}\tau}\frac{\sinh(n_{\lambda}+1)\tau}{\sinh \tau}=
\sum_{k=0}^{n_{\lambda}}x^{k},\quad x=e^{-2\tau},\\
f(\tau)=
\frac{(n_{\lambda}+1)\sinh (\tau) \cosh ((n_{\lambda}+1)\tau)-\cosh (\tau)\sinh ((n_{\lambda}+1)\tau)}
{2 \sinh (\tau) \sinh ((n_{\lambda}+1)\tau)}
-\frac{n_{\lambda}}{2}.
\end{gather*}
Furthermore, for $0 \leq \tau$ if and only if $-\frac{n_{\lambda}}{2} \leq f(\tau) < 0$,
and $f(0)=-\frac{n_{\lambda}}{2}$.
\label{nonnegu2}
\end{lem}

\begin{proof}
Since we have $h'(\tau)=\ispa{(\partial k)(\tau \alpha),\alpha}$ and $\ispa{\alpha,\alpha}=2$,
the differential $(\partial k)(\tau \alpha)$ is given by
$(\partial k)(\tau \alpha)=\frac{1}{2}h'(\tau)\alpha$. The equation
\eqref{momentu2} follows from this and the definition of the moment map.
The explicit expression for the function $h(\tau)$ follows from \eqref{schurch} and \eqref{schuru2},
and that for $f(\tau)$ is shown by a direct computation.
Next, it is easy to show that, by using the expression for $h(\tau)$ in terms of a
polynomial in $x=e^{-2\tau}$, $f(0)=n_{\lambda}/2$.
Also, we have $\lim_{\tau \to +\infty}f(\tau)=0$ and $\lim_{\tau \to -\infty}f(\tau)=n_{\lambda}$.
From this the rest of the assertion follows.
\end{proof}

Finally, we shall examine that the term of the Weyl denominator in Theorem \ref{IMULT} does not
vanish for generic dominant weight in the case where $G=U(2)$.

\begin{prop}
Let $\nu_{j}$ $(1 \leq j \leq n_{\lambda}/2)$ be a dominant weight defined in \eqref{weightu2}.
We set $\tau_{j}:=\tau_{\lambda}(\nu_{j})$: $\tau_{j}$ is the unique non-negative number satisfying
$f(\tau_{j})=-j$, where $f(\tau)$ is defined by \eqref{momentu2}.
Then the multiplicity $a_{N}(\lambda;N\nu_{j})$ of $V_{N\nu_{j}}$ in $V_{\lambda}^{\otimes N}$ has
the following asymptotic formula:
\[
a_{N}(\lambda;N\nu)=(2\pi N)^{-1/2}
e^{-N(n_{\lambda}-2j)}
\left(
\frac{\sinh (n_{\lambda}+1)\tau_{j}}{\sinh \tau_{j}}
\right)^{N}
\left(
a_{\lambda}(j) +O(N^{-1})
\right),
\]
where the positive constant $a_{\lambda}(j)$ is given by
\[
a_{\lambda}(j)=2e^{-\tau_{j}}
\sqrt{
\frac{2\sinh^{4}\tau_{j}\sinh^{2}(n_{\lambda}+1)\tau_{j}}
{\sinh^{2}(n_{\lambda}+1)\tau_{j}-(n_{\lambda}+1)^{2}\sinh^{2}\tau_{j}}
}.
\]
The leading term $a_{j}$ vanishes if and only if
$n_{\lambda}$ is even and $j=n_{\lambda}/2$. In this case, the dominant weight $\nu_{j}$ $(j=n_{\lambda}/2)$
is in the unique wall of the Weyl chamber $C$.
\end{prop}

\begin{proof}
The non-negativity of the number $\tau_{j}$ follows form Lemma \ref{nonnegu2} and
that $\nu_{j}$ is a dominant weight, {\it i.e.}, $1 \leq j \leq n_{\lambda}/2$.
The lattice $L^{*}=X^{*} \cap I^{*} =X^{*} \cap \mb{Z}^{2}$ is spanned by the simple
root $\alpha$. Thus we have $\Lambda =L$, and hence the finite group $\Pi(U(2))$ is trivial.
Note that the Weyl denominator $\Delta (\tau \alpha/2\pi i)$ is given by
\[
\Delta (\tau \alpha /2\pi i)=2\sinh\tau,
\]
which is non-negative for $\tau=\tau_{j}$ and zero if and only if $\tau=0=\tau_{n_{\lambda}/2}$.
By \eqref{deltau2}, the positive constant $\delta_{\lambda}(\nu_{j})$ is given by
\[
e^{\delta_{\lambda}(\nu_{j})}=h(\tau_{j})^{N}e^{2j\tau_{j}}
=e^{-(n_{\lambda}-2j)\tau_{j}}
\left(
\frac{\sinh (n_{\lambda}+1)\tau_{j}}{\sinh\tau_{j}}
\right).
\]
Note that half the sum of the positive roots is given by $\rho =\alpha/2$, and hence
$\ispa{\rho,\tau_{j}\alpha}=\tau_{j}$.
Recall that the matrix $A_{\lambda}(\nu_{j})$ is equal to $A(\tau_{\lambda}(\nu))$ where
$A(\tau)$ ($\tau \in X$) is the derivative of the moment map $\mu_{P}(\tau)$.
In our case, $A(\tau)$ is a positive real number given by
\[
A(\tau) = \frac{h(\tau)h''(\tau) -h'(\tau)^{2}}
{2h(\tau)^{2}}
=\frac{\sinh^{2}(n_{\lambda}+1)\tau-(n_{\lambda}+1)^{2}\sinh^{2}\tau}
{2\sinh^{2}\tau \sinh^{2}(n_{\lambda}+1)\tau}.
\]
(Note that, since $\ispa{\alpha,\alpha}=2$,
$\alpha \otimes \alpha$ is identified with the multiplication by $2$.)
Therefore, the assertion follows from Theorem \ref{IMULT}.
\end{proof}

\section{\label{FINAL} Final Comments}

We close with some remarks on lattice paths and also on the
symplectic interpretation of our problems and results.

\subsection{Further relations between multiplicities of
irreducibles and lattice paths}

A number of relations are known between   lattice path counting
problems to that of determining multiplicities of weights in
tensor powers $V_{\lambda}^{\otimes N}$. We used the formulae
(\ref{MULTM}) and (\ref{MULTA}) in terms of weighted
multiplicities of lattice paths. There are other formulae which
express multiplicities in terms of unweighted but constrained
sums.

One  is given by Theorem 2 of the paper of Grabiner-Magyar
\cite{GM}: {\it Let  $C$ be the Weyl chamber of a reductive
complex Lie algebra, $V$ be a finite dimensional representation,
$S$ be the set of weights of $V$ and $L$ be a lattice containing
$S$ and $\rho$. Then the number $b_{\rho, \rho + \mu, N}$ of walks
of $N$ steps from $\rho$ to $\rho + \mu$ which stay strictly
within $C$ equals the multiplicity of the irreducible with highest
weight $\mu$ in $V^{\otimes N}$.}

To use this formula, one needs to count lattice paths satisfying
the constraint.  One possible approach is to use the
Gessel-Zeilberger formula, which
  says that for
{\it reflectable} paths, one has:
\[\label{GZ}
b_{\eta , \lambda, N} = \sum_{w \in W} sgn(w)
\mathcal{P}_{N}(\lambda -w\eta),.
\]
This provideds an alternative approach to our problems, but
requires an analysis of reflectable paths.
 Many further (and much more general) relations between characters and multiplicities
to sums over special lattice paths are discussed in \cite{Lit}.

\subsection{\label{SYMPMOD} Symplectic model}

The reader may note a  resemblance between the problems studied in
this paper and   the well-known problem of finding asymptotics of
weight
 multiplicities in
$V_{N \lambda}$, where $V_{N \lambda}$ is the irreducible with
highest weight $N \lambda$ (see e.g. \cite{H, GS}). In both cases,
the possible weights lie in $Q(N \lambda)$ and one may define
analogous distribution of weights of $V_{N \lambda}$. However, the
relation is not very close, since our problem is about the
thermodynamic limit rather than the semiclassical limit. We add a
few remarks to clarify the relations.

 We recall
the symplectic interpretation of the latter multiplicity problem:
the maximal torus ${\bf T}$ acts by conjugation on the co-adjoint
orbit ${\mathcal O}_{\lambda}$ associated to $V_{\lambda}$  in a
Hamiltonian fashion,  with moment map given by the orthogonal
projection $\mu_{\lambda}: {\mathcal O}_{\lambda} \to {\bf t}^*$
to the Cartan dual subalgebra. The image is given by
$\mu_{\lambda} ({\mathcal O}_{\lambda}) = Q(\lambda)$. As proved
by G. Heckman, multiplicities of weights in $V_{N \lambda}$ become
asymptotically distributed according to the (Duistermaat-Heckman)
measure, namely the push-forward $\mu_{\lambda *} dVol_{\lambda}$
the symplectic volume measure of ${\mathcal O}_{\lambda}$ under
the orthogonal projection to ${\bf t}^*$ \cite{H, GS}.

The limit formula in Theorem  \ref{SIMPLE} also  has a symplectic
interpretation: To $V_{\lambda}^{\otimes N}$ corresponds the
symplectic manifold
$${\mathcal O}_{\lambda}^N:= {\mathcal O}_{\lambda} \times \cdots
\times {\mathcal O}_{\lambda}\;\;\; (N \; \mbox{times}).$$ Then
${\bf T}$ acts on ${\mathcal O}_{\lambda}^N$ with moment map
\begin{equation} \mu_{\lambda}^N: {\mathcal O}_{\lambda}^N \to
{\bf t}^*,\;\;\;\mu_{\lambda}^N (x_1, \dots, x_N) = \mu_{\lambda}
(x_1) + \cdots +  \mu_{\lambda}(x_N). \end{equation} The image of
the moment map is the convex polytope $Q(N \lambda) = N
\mu({\mathcal O}_{\lambda})$, and one may define  the
Duistermaat-Heckman type measure on $Q(\lambda)$ by:
\begin{equation} \label{CLASSICAL}
dm_{\lambda}^N := D_{N}^{-1}  (\mu_{\lambda}^{N})_*
(dVol_{\lambda} \times \cdots \times dVol_{\lambda} )\;\;\; (N \;
\mbox{times}) \end{equation} on $Q(\lambda)$, where $D_N x = N x$
is the dilation operator. Equivalently, this latter measure is
defined by
\begin{equation} \int_{Q(\lambda)} f(x) dm_{\lambda}^N(x) = \int_{{\mathcal O}_{\lambda} \times \cdots
\times {\mathcal O}_{\lambda}}  f(\frac{\mu_{\lambda} (x_1) +
\cdots + \mu_{\lambda}(x_N)}{N})
 dVol_{\lambda}(x_1) \times \cdots \times dVol_{\lambda}(x_N).
 \end{equation}
 Thus, $dm_{\lambda}^N$ is the distribution of the sum of the  (vector valued)
 independent random variables $\mu_{\lambda}(x_j)$, the law of large numbers implies  that the limit
 equals the mean value of the random variables:
\begin{equation} \label{LLN} dm_{\lambda}^N \; \to
\delta_{Q^{*}(\lambda)},\;\; \mbox{weakly as}\;\; N \to \infty.
\end{equation}

This measure represents the thermodynamic limit of the classical
spin chain with phase space ${\mathcal O}_{\lambda}$ at each site,
while our problem involves the thermodynamic limit of the quantum
spin chain. The two problems are quite distinct until one lets the
weight $\lambda \to \infty$ along a ray, i.e. considers the joint
asymptotics of weights in  $V_{M \lambda}^{\otimes N}$. The
Heckman theorem says that if $N$ is fixed and $M \to \infty$ then
the quantum problem converges to the classical one. It would be
interesting to investigate the joint asymptotics as both
parameters become large.

\end{document}